\documentclass[11pt]{amsart}

\usepackage[all]{xy}
\usepackage{mathrsfs}
\usepackage{amsmath,amsthm, amssymb, amsgen,amsxtra,amsfonts, amsbsy} 
\usepackage[all]{xy}
\usepackage{verbatim}
\usepackage{array}
\usepackage{hyperref}
\usepackage{graphicx}
\usepackage{caption,rotating}
\usepackage{color}
\usepackage{subcaption}
\usepackage[shortlabels]{enumitem}
\usepackage{times}
\newcommand{\mathds}[1]{{\mathbb #1}}
\usepackage{stmaryrd}
\usepackage[graphicx]{realboxes}

\begin{document}
%
%
%
\theoremstyle{definition}
\newtheorem{Definition}{Definition}[section]
\newtheorem*{Definitionx}{Definition}
\newtheorem*{WeierstrassModelx}{Weierstra{\ss} Model}
\newtheorem*{Conventionx}{Convention}
\newtheorem{Construction}{Construction}[section]
\newtheorem{Example}[Definition]{Example}
\newtheorem{Examples}[Definition]{Examples}
\newtheorem{Remark}[Definition]{Remark}
\newtheorem*{Remarkx}{Remark}
\newtheorem{Remarks}[Definition]{Remarks}
\newtheorem{Caution}[Definition]{Caution}
\newtheorem{Conjecture}[Definition]{Conjecture}
\newtheorem*{Conjecturex}{Conjecture}
\newtheorem{Question}[Definition]{Question}
\newtheorem*{Questionx}{Question}
\newtheorem*{Acknowledgements}{Acknowledgements}
\newtheorem*{Organization}{Organization}
\newtheorem*{Disclaimer}{Disclaimer}
\newtheorem{Caveat}[Definition]{Caveat}
\theoremstyle{plain}
\newtheorem{Theorem}[Definition]{Theorem}
\newtheorem*{Theoremx}{Theorem}
\newtheorem{innercustomthm}{Theorem}
\newenvironment{customthm}[1]
  {\renewcommand\theinnercustomthm{#1}\innercustomthm}
  {\endinnercustomthm}
\newtheorem{Proposition}[Definition]{Proposition}
\newtheorem*{Propositionx}{Proposition}
\newtheorem{Lemma}[Definition]{Lemma}
\newtheorem{Corollary}[Definition]{Corollary}
\newtheorem*{Corollaryx}{Corollary}
\newtheorem{Fact}[Definition]{Fact}
\newtheorem{Facts}[Definition]{Facts}
\newtheoremstyle{voiditstyle}{3pt}{3pt}{\itshape}{\parindent}%
{\bfseries}{.}{ }{\thmnote{#3}}%
\theoremstyle{voiditstyle}
\newtheorem*{VoidItalic}{}
\newtheoremstyle{voidromstyle}{3pt}{3pt}{\rm}{\parindent}%
{\bfseries}{.}{ }{\thmnote{#3}}%
\theoremstyle{voidromstyle}
\newtheorem*{VoidRoman}{}

%
\newcommand{\prf}{\par\noindent{\sc Proof.}\quad}
\newcommand{\blowup}{\rule[-3mm]{0mm}{0mm}}
\newcommand{\cal}{\mathcal}
\newcommand{\Aff}{{\mathds{A}}}
\newcommand{\BB}{{\mathds{B}}}
\newcommand{\CC}{{\mathds{C}}}
\newcommand{\EE}{{\mathds{E}}}
\newcommand{\FF}{{\mathds{F}}}
\newcommand{\GG}{{\mathds{G}}}
\newcommand{\HH}{{\mathds{H}}}
\newcommand{\NN}{{\mathds{N}}}
\newcommand{\ZZ}{{\mathds{Z}}}
\newcommand{\PP}{{\mathds{P}}}
\newcommand{\QQ}{{\mathds{Q}}}
\newcommand{\RR}{{\mathds{R}}}
\newcommand{\Liea}{{\mathfrak a}}
\newcommand{\Lieb}{{\mathfrak b}}
\newcommand{\Lieg}{{\mathfrak g}}
\newcommand{\Liem}{{\mathfrak m}}
\newcommand{\ideala}{{\mathfrak a}}
\newcommand{\idealb}{{\mathfrak b}}
\newcommand{\idealg}{{\mathfrak g}}
\newcommand{\idealm}{{\mathfrak m}}
\newcommand{\idealp}{{\mathfrak p}}
\newcommand{\idealq}{{\mathfrak q}}
\newcommand{\idealI}{{\cal I}}
\newcommand{\lin}{\sim}
\newcommand{\num}{\equiv}
\newcommand{\dual}{\ast}
\newcommand{\iso}{\cong}
\newcommand{\homeo}{\approx}
\newcommand{\mm}{{\mathfrak m}}
\newcommand{\pp}{{\mathfrak p}}
\newcommand{\qq}{{\mathfrak q}}
\newcommand{\rr}{{\mathfrak r}}
\newcommand{\pP}{{\mathfrak P}}
\newcommand{\qQ}{{\mathfrak Q}}
\newcommand{\rR}{{\mathfrak R}}
%
%
\newcommand{\OO}{{\cal O}}
\newcommand{\numero}{{n$^{\rm o}\:$}}
\newcommand{\mf}[1]{\mathfrak{#1}}
\newcommand{\mc}[1]{\mathcal{#1}}
\newcommand{\into}{{\hookrightarrow}}
\newcommand{\onto}{{\twoheadrightarrow}}
\newcommand{\Spec}{{\rm Spec}}
\newcommand{\BigSpec}{{\rm\bf Spec}\:}
\newcommand{\Spf}{{\rm Spf}\:}
\newcommand{\Proj}{{\rm Proj}\:}
\newcommand{\Pic}{{\rm Pic }}
\newcommand{\MW}{{\rm MW }}
\newcommand{\Br}{{\rm Br}}
\newcommand{\NS}{{\rm NS}}
\newcommand{\CH}{{\rm CH}}
\newcommand{\Sym}{{\mathfrak S}}
\newcommand{\Lie}{{\rm Lie}}
\newcommand{\Aut}{{\rm Aut}}
\newcommand{\Num}{{\rm Num}}
\newcommand{\Autp}{{\rm Aut}^p}
\newcommand{\Hom}{{\rm Hom}}
\newcommand{\Ext}{{\rm Ext}}
\newcommand{\ord}{{\rm ord}}
\newcommand{\Char}{{\rm char}}
\newcommand{\rank}{{\rm rank}}
\newcommand{\coker}{{\rm coker}\,}

\newcommand{\sm}{{\rm sm}}
\newcommand{\divisor}{{\rm div}}
\newcommand{\Def}{{\rm Def}}
\newcommand{\piet}{{\pi_1^{\rm \acute{e}t}}}
\newcommand{\Het}[1]{{H_{\rm \acute{e}t}^{{#1}}}}
\newcommand{\Hfl}[1]{{H_{\rm fl}^{{#1}}}}
\newcommand{\Hcris}[1]{{H_{\rm cris}^{{#1}}}}
\newcommand{\HdR}[1]{{H_{\rm dR}^{{#1}}}}
\newcommand{\hdR}[1]{{h_{\rm dR}^{{#1}}}}

\newcommand{\loc}{{\rm loc}}
\newcommand{\et}{{\rm \acute{e}t}}
\newcommand{\I}{{\rm I}}
\newcommand{\II}{{\rm II}}
\newcommand{\defin}[1]{{\bf #1}}
\newcommand{\oX}{\cal{X}}
\newcommand{\oA}{\cal{A}}
\newcommand{\oY}{\cal{Y}}
\newcommand{\blue}{\textcolor{blue}}
\newcommand{\red}{\textcolor{red}}

\title{Automorphisms of unnodal Enriques surfaces}

\author{Gebhard Martin}
\address{\hfill \newline
Mathematisches Institut  \newline
Universit\"at Bonn \newline
Endenicher Allee 60 \newline
53115 Bonn \newline
Germany}
\email{gmartin@math.uni-bonn.de}

\date{\today}

\begin{abstract}
It follows from an observation of A. Coble in 1919 that the automorphism group of an unnodal Enriques surface contains the $2$-congruence subgroup of the Weyl group of the $E_{10}$-lattice. In this article, we determine how much bigger the automorphism group of an unnodal Enriques surface can be. Furthermore, we show that the automorphism group is in fact equal to the $2$-congruence subgroup for generic Enriques surfaces in arbitrary characteristic (under the additional assumption that the Enriques surface is ordinary if the characteristic is $2$), improving the corresponding result of W. Barth and C. Peters for very general Enriques surfaces over the complex numbers.
\end{abstract}
\maketitle

\begin{Conventionx}
Except for Section \ref{genusone}, we will be working over an algebraically closed field $k$. If $\Char(k) = 2$, all Enriques surfaces considered will be assumed to be ordinary. In particular, their canonical cover is a smooth K3 surface. 
\end{Conventionx}

Denote by $W_{E_{10}}$ the Weyl group of $E_{10} = U \oplus E_8$, which is the unique even unimodular lattice of signature $(1,9)$. It follows from the classification of Enriques surfaces with numerically trivial automorphisms (see \cite{MukaiNamikawa} and \cite{DolgachevMartin}), that $\Aut(X)$ is a subgroup of $W_{E_{10}}$ for every \emph{unnodal} Enriques surface $X$, that is, for every Enriques surface which does not contain $(-2)$-curves. On the other hand,  an observation of A. Coble \cite{Coble} (see \cite{Allcock} for an elegant proof in arbitrary characteristic) implies that for an unnodal Enriques surface $X$ the group $\Aut(X)$ always contains the (infinite) $2$-congruence subgroup of $W_{E_{10}}$,
$$
W_{E_{10}}(2) := \ker(W_{E_{10}} \to {\rm O}(E_{10}/2E_{10})),
$$
as a normal subgroup. Hence, we can define the \emph{extra automorphism group} $$\overline{\Aut}(X) :=\Aut(X)/(W_{E_{10}}(2))$$ of an unnodal Enriques surface $X$. The purpose of this article is to study the group $\overline{\Aut}(X)$. The following are our main results.
%
%
%

\begin{customthm}{A}\label{mainA}
Let $X$ be an unnodal Enriques surface. Then, $|\overline{\Aut}(X)| \in \{1,2,4\}$. If $\Char(k) = 2$, even $|\overline{\Aut}(X)| = \{1\}$ holds.
\end{customthm}

\begin{customthm}{B}\label{mainB}
Let $X$ be a generic Enriques surface. Then, $X$ is unnodal and $\overline{\Aut}(X)$ is trivial, that is, $\Aut(X) = W_{E_{10}}(2)$. 
\end{customthm}

More precisely, we show that the locus of Enriques surfaces with extra automorphisms is contained in a closed subset of codimension $5$ in the base of any modular family of unnodal Enriques surfaces (see Corollary \ref{CorollaryB} for a more refined statement).

\begin{Remark}
All the possibilities for $|\overline{\Aut}(X)|$ given in Theorem \ref{mainA} are actually realized on unnodal Enriques surfaces. More precisely, in Section \ref{section6}, we will give examples, first described by H. Ohashi \cite{Ohashi} over the complex numbers in a different context, realizing all the possibilities in Theorem \ref{mainA} in all except possibly finitely many characteristics. These families of examples are $5$- resp. $2$-dimensional, showing that the bounds on the dimension of such families given in Corollary \ref{CorollaryB} are sharp.
\end{Remark}

\begin{Remark} We note that, even for $k = \CC$, Theorem \ref{mainA} is new. Moreover, over the complex numbers, Theorem \ref{mainB} extends the main result of W. Barth and C. Peters in \cite{BarthPeters} (which also follows from independent work of V. V. Nikulin \cite{Nikulin}) according to which $|\overline{\Aut}(X)|= 1$ holds for a very general Enriques surface, whereas we show that this in fact holds outside a closed subset of codimension $5$ in the moduli space of unnodal Enriques surfaces.
\end{Remark}

\begin{Remark}\label{Remark2}
In analogy to the introduction of \cite{BarthPeters}, we deduce that there are exactly $527 = 17 \cdot 31$ elliptic fibrations on a generic Enriques surface $X$ up to automorphisms of $X$. As in \cite{BarthPeters}, it also follows that there are exactly $2^7 \cdot 17 \cdot 31$ double plane models, $2^{11} \cdot 5 \cdot 17 \cdot 31$ Enriques models, and $2^{14} \cdot 3 \cdot 17 \cdot 31$ Fano models of $X$ up to the action of $\Aut(X)$ for a generic Enriques surface $X$ in arbitrary characteristic.
\end{Remark}

The outline of this paper is as follows: After recalling the known facts about automorphisms of unnodal Enriques surfaces in Section \ref{section1}, we study the action of the automorphism group of an elliptic fibration on its Jacobian fibration in Section \ref{section2}. Due to the level of generality of the exposition, this may be of independent interest. The Jacobian fibration of an elliptic fibration of an Enriques surface is a rational elliptic surface and therefore we gather all the facts about rational elliptic surfaces that we need in Section \ref{section3}. Finally, in Sections \ref{section4} and \ref{section5}, we apply the machinery of Section \ref{section2} together with the facts collected in Section \ref{section3} to unnodal Enriques surfaces and deduce Theorem \ref{mainA} and Theorem \ref{mainB}. After that, in Section \ref{section6}, we will give examples of unnodal Enriques surfaces with extra automorphisms, realizing all the possibilities for $|\overline{\Aut}(X)|$ given in Theorem \ref{mainA}.

 \medskip
\noindent
{\bf Acknowledgement.}
This project started as a joint project with Shigeyuki Kond\={o} and I am very grateful to him for many stimulating discussions. I would like to thank Claudia Stadlmayr for numerous important remarks and Simon Brandhorst for pointing me to the work of Junmyeong Jang. Finally, I would also like to thank Christian Liedtke and Matthias Sch\"utt for helpful comments on a first version of this article.

\section{Generalities on unnodal Enriques surfaces} \label{section1}
In this section, we will collect general facts about unnodal Enriques surfaces and their automorphism groups which we will need later on. Throughout this section, $X$ will denote an unnodal Enriques surface.

\subsection{Elliptic fibrations on unnodal Enriques surfaces} \label{EllfibronEnriques}
Recall the following facts on elliptic fibrations $f: X \to C$ of an unnodal Enriques surface $X$ (e.g. from \cite[Chapter V]{CossecDolgachev} and \cite{LiuLorenziniRaynaud}):
\begin{itemize}
\item $C \cong \PP^1$.
\item All fibers of $f$ are irreducible. In particular, the fibers of $f$ are of type $\I_0$, $\I_1$ or $\II$ in Kodaira's notation for the fibers of elliptic fibrations.
\item $f$ has exactly two double fibers if $\Char(k) \neq 2$ (resp. exactly one if $\Char(k) = 2$). A reduced curve with support equal to the support of a double fiber is called a \emph{half-fiber} of $f$.
\item $f$ admits a separable bisection, that is, there is an irreducible curve on $X$ meeting a general fiber of $f$ in two distinct points.
\item The half-fibers of $f$ are of type $\I_0$ or $\I_1$ (if $\Char(k) = 2$, the half-fiber of $f$ is of type $\I_0$ and ordinary, or of type $\I_1$).
\item The Jacobian $J(f): J \to \PP^1$ of $f$ is a rational elliptic surface such that every fiber of $J(f)$ is of the same type as the corresponding fiber of $f$.
\end{itemize}

\begin{Remark}\label{atleastonesingfiber}
Note that, since the discriminant of $J(f)$ is a polynomial of degree $12$ which vanishes exactly once at a fiber of type $\I_1$ and does not vanish at a fiber of type $\I_0$, $f$ and $J(f)$ have either at least one fiber of type $\II$ or twelve fibers of type $\I_1$. In particular, there is at least one singular fiber of $f$ which is not a double fiber.
\end{Remark}

\subsection{Automorphisms of unnodal Enriques surfaces}
Recall the following facts and definitions related to the numerical lattice $\Num(X)$ of $X$ (from \cite{BarthPeters} and \cite{CossecDolgachev}):
\begin{itemize}
\item $\Num(X) = E_{10} := U \oplus E_8$, where $U$ is the hyperbolic plane and $E_8$ is the unique even unimodular negative definite lattice of rank $8$. We will fix one such identification for the remainder of this article.
\item The Weyl group $W_{E_{10}}$ of $E_{10}$ is the group of isometries of $E_{10}$ preserving the two half-cones of the positive cone.
\item The $2$-congruence subgroup $W_{E_{10}}(2)$ of $W_{E_{10}}$ is defined as
$$
W_{E_{10}}(2) := \ker(W_{E_{10}} \to {\rm O}(E_{10}/2E_{10})).
$$
\item We have $W_{E_{10}}/(W_{E_{10}}(2)) \cong {\rm O}^+(10,\mathbb{F}_2)$.
\end{itemize}

Since the action of $\Aut(X)$ on $\Num(X)$ preserves the half-cones of the positive cone and $ W_{E_{10}} \subseteq {\rm O}(E_{10})$ is the group of isometries preserving these half-cones, there is a homomorphism $\Aut(X) \to W_{E_{10}}$. The kernel of this map has been studied by S. Mukai and Y. Namikawa \cite{MukaiNamikawa} over the complex numbers and by I. Dolgachev and the author \cite{DolgachevMartin} in positive characteristic. The relevant result in our situation is the following.

\begin{Theorem} [\cite{MukaiNamikawa}, \cite{DolgachevMartin}] \label{numtriv}
An unnodal Enriques surface $X$ has no numerically trivial automorphisms, that is, the group $\Aut_{nt}(X) := \ker(\Aut(X) \to W_{E_{10}})$ is trivial. In particular, we can identify $\Aut(X)$ with a subgroup $\Aut(X) \subseteq W_{E_{10}}$ of $W_{E_{10}}$.
\end{Theorem}

The simplest automorphisms of unnodal Enriques surfaces arise as follows.

\begin{Example}\label{BiellipticInvolution}
A bielliptic (or Bertini-type) involution of an unnodal Enriques surface $X$ is the covering involution $g$ associated to a linear system $|2F_1 + 2F_2|$, where $F_1$ and $F_2$ are half-fibers of elliptic fibrations of $X$ with $F_1.F_2 = 1$. The map induced by $|2F_1 + 2F_2|$ is a finite morphism of degree $2$ onto a quartic symmetroid del Pezzo surface (see \cite[Section IV.7]{CossecDolgachev} or \cite{DolgachevMartin}). Note that the covering involution $g$ acts trivially on the pencils $|2F_1|$ and $|2F_2|$. Moreover, the pair $(F_1,F_2)$ yields an embedding $\langle F_1,F_2 \rangle \cong U \hookrightarrow E_{10}$ with orthogonal complement $E_8$. In the induced decomposition $E_{10} \cong U \oplus E_8$, $g$ acts as ${\rm id}_U$ on $U$ and as $-{\rm id}_{E_8}$ on $E_8$.
\end{Example}

In particular, we see that every bielliptic involution is contained in $W_{E_{10}}(2)$.
In fact, for unnodal Enriques surfaces $X$, an observation of A. Coble shows that $W_{E_{10}}(2)$ is generated by bielliptic involutions.

\begin{Theorem}[\cite{Coble}, \cite{Allcock}] \label{Allcock}
Let $X$ be an unnodal Enriques surface. Then, $W_{E_{10}}(2) \subseteq \Aut(X)$ is the normal subgroup generated by bielliptic involutions.
\end{Theorem}

From Theorem \ref{numtriv} and Theorem \ref{Allcock} we deduce that $W_{E_{10}}(2) \subseteq \Aut(X) \subseteq W_{E_{10}}$.
Furthermore, we have $W_{E_{10}}/W_{E_{10}}(2) \cong {\rm O}^+(10,\mathbb{F}_2)$. The order of this group is $|{\rm O}^+(10,\mathbb{F}_2)| = 2^{21} \cdot 3^5 \cdot 5^2 \cdot 7 \cdot 17 \cdot 31$. For future reference, let us record the consequences of this observation for the possible orders of $\overline{\Aut}(X)$ as a corollary.

\begin{Corollary}\label{possibleorders}
For every unnodal Enriques surface $X$, the group $\overline{\Aut}(X)$ is a subgroup of ${\rm O}^+(10,\mathbb{F}_2)$. In particular, its order divides $2^{21} \cdot 3^5 \cdot 5^2 \cdot 7 \cdot 17 \cdot 31$.
\end{Corollary}

\begin{Caveat}
The group ${\rm O}^+(10,\mathbb{F}_2)$ is the automorphism group of a non-degenerate symmetric bilinear form on $\mathbb{F}_2^{10}$ of Witt defect $0$ (i.e., containing a totally isotropic subspace of dimension $5$). It contains a simple group of index $2$ which is sometimes referred to using the same notation (e.g. in the Atlas \cite{Atlas}).
\end{Caveat}

In \cite{BarthPeters}, W. Barth and C. Peters studied the action of $W_{E_{10}}$ on $E_{10}$ to compute the automorphism group of a very general complex Enriques surface. Moreover, they studied geometric realizations of very general complex Enriques surfaces by computing the stabilizers of the $W_{E_{10}}$-action on sets of isotropic sequences of the following form.

\begin{Definition}
A (\emph{non-degenerate canonical}) $U_{[n]}$-\emph{sequence} on $X$ is a sequence $(f_1,\hdots,f_n)$ with $f_i \in \Num(X)$ such that
\begin{itemize}
\item $f_i.f_j = 1 - \delta_{ij}$,
\item the $f_i$ are nef and primitive.
\end{itemize}
\end{Definition}

\begin{Example}\label{examples}
Let $(F_1,\hdots,F_n)$ be a lift to $\Pic(X)$ of a $U_{[n]}$-sequence $(f_1,\hdots,f_n)$ on an unnodal Enriques surface $X$. Then, the sequence $(F_1,\hdots,F_n)$ gives rise to some of the most classical geometric models of $X$ as follows (see e.g. \cite[Chapter IV]{CossecDolgachev} and \cite[Section 2]{DolgachevMark}):

\begin{itemize}
\item Elliptic fibrations:  $|2F_1|$ induces an elliptic fibration of $X$ with half-fiber $F_1$ (and $F_1' \sim F_1 + K_X$ if $\Char(k) \neq 2$).
\smallskip
\item Double plane models: 
 $|2F_1 + 2F_2|$ induces a morphism of degree $2$ from $X$ onto a quartic symmetroid del Pezzo as in Example \ref{BiellipticInvolution}.
 \smallskip
\item Enriques models: 
$|F_1 + F_2 + F_3|$ induces a birational map from $X$ onto an Enriques sextic in $\PP^3$.
\smallskip
\item Fano models: 
 $|\frac{1}{3}\sum_{i =1}^{10} F_i|$ induces an isomorphism of $X$ with an intersection of $10$ cubics in $\PP^5$.
\end{itemize}
\end{Example}

While the group $W_{E_{10}}$ acts transitively on the set of all $U_{[n]}$-sequences for $n =1,2,3,10$, this is no longer true for $W_{E_{10}}(2)$. However, the number of $W_{E_{10}}(2)$-orbits of $U_{[n]}$-sequences has been calculated explicitly in these cases in \cite[p. 397]{BarthPeters}. We remark that the following proposition is a statement about the action of the group $W_{E_{10}}(2)$ on the lattice $E_{10}$ and in particular it does not depend on the characteristic of $k$.

\begin{Proposition}[\cite{BarthPeters}] \label{orbits}
The number $N(n)$ of $W_{E_{10}}(2)$-orbits of elements of the form $\sum_{i=1}^n f_i$, where $(f_1,\hdots,f_n)$ is a non-degenerate canonical $U_{[n]}$-sequence, is as given in the following table.

\begin{table}[h]
$\begin{array}{|l|l|} \hline
n & N(n) \\ \hline \hline
1 & 17 \cdot 31 = 527 \\
2 & 2^7 \cdot 17 \cdot 31 = 67456 \\
3 & 2^{10}\cdot 5\cdot 17 \cdot 31 = 2698240 \\
10 & 2^{13} \cdot 3 \cdot 17 \cdot 31 = 12951552 \\ \hline
\end{array}$
\vspace{1mm}
\caption{Number $N(n)$ of $W_{E_{10}}(2)$-orbits of sums over $U_{[n]}$-sequences}
\label{table1}
\end{table}
\end{Proposition}

\vspace{-8mm}
\subsection{Lattice theoretic preliminaries}
To obtain restrictions on $\overline{\Aut}(X)$, we will make use of the K3 cover $\pi: Y \to X$ of $X$ at several points in this article. An isometry $g \in {\rm O}(\Num(X))$ induces an isometry $\pi^* g := \frac{1}{2}\pi^* \circ g \circ \pi_*$ of $\pi^*\Num(X) \cong E_{10}(2)$. For the following lemma, see \cite[Section 1.2]{BarthPeters}.

\begin{Lemma} \label{identification}
Let $X$ be an unnodal Enriques surface and let $\pi: Y \to X$ be its K3 cover. The natural map $\pi^*: {\rm O}(\Num(X)) \to {\rm O}(\pi^*(\Num(X))) = {\rm O}(E_{10}(2))$
induces an isomorphism $W_{E_{10}}/(W_{E_{10}}(2)) \cong {\rm O}(E_{10}(2)^{\vee}/E_{10}(2))$. In particular, $\overline{\Aut}(X)$ acts faithfully on the discriminant group of $\pi^*(\Num(X))$.
\end{Lemma}
%
%

Therefore, every automorphism $g \in \Aut(X)$ which is non-trivial modulo $W_{E_{10}}(2)$ acts non-trivially on $E_{10}(2)^{\vee}/E_{10}(2)$. On the other hand, $g$ is induced by an automorphism of $Y$, hence $g^* \in {\rm O}(\pi^*(\Num(X)))$ extends to an isometry of $\Num(Y)$ and in fact to every cohomology group of $Y$ containing $\pi^*(\Num(X))$. The following criterion for extendability of lattice automorphisms shows that we can use this observation to obtain restrictions on $g$ and thus on $\overline{\Aut}(X)$.

\begin{Lemma} \label{extendability} \cite[Corollary 1.5.2.]{Nikulin}
Let $M$ be an even lattice and let $S \subseteq M$ be a primitive sublattice with orthogonal complement $K$. Let $H = M/(S \oplus K) \subseteq (S \oplus K)^\vee / (S \oplus K)$. The group $H$ is identified with a subgroup $H_S$ of $S^\vee/S$ and $H_K$ of $K^\vee/K$ via the two projections. We denote by $\gamma: H_S \xrightarrow{\sim} H_K$ the induced isomorphism.

Let $\alpha: S \to S$ and $\beta: K \to K$ be isometries. The isometry $(\alpha,\beta)$ of $S \oplus K$ lifts to $M$ if and only if $\alpha$ preserves $H_S$, $\beta$ preserves $H_K$ and the induced automorphisms $\tilde{\alpha}: H_S \to H_S$ and $\tilde{\beta}:H_K \to H_K$ satisfy 
\begin{equation*}
\tilde{\beta} \circ \gamma = \gamma \circ \tilde{\alpha}.
\end{equation*}
\end{Lemma}

\begin{Remark}
If $M$ is unimodular, we have $H_S = S^\vee/S$ and $H_K = K^\vee/K$, so the first two conditions of Lemma \ref{extendability} are automatically satisfied.
\end{Remark}


%


\section{Automorphisms of elliptic fibrations} \label{section2}  \label{genusone}
In this section, we will develop some general facts on automorphisms of curves of genus one over a not necessarily algebraically closed field $K$ in order to be able to apply them to elliptic surfaces over algebraically closed fields. In particular, we will study how automorphisms of a curve of genus one act on its Picard scheme.
The results of this section will be applied to elliptic fibrations of unnodal Enriques surfaces.
\subsection{Automorphisms of curves of genus one}
We fix the following notations and identifications:

\begin{itemize}
\item $C = \Pic^1_{C/K}$ is a smooth curve of genus one over a field $K$. 
\item $(E,O) := (\Pic^0_{C/K},\mathcal{O}_C)$ is its Jacobian.
\item $\Aut_{C/K}$ is the automorphism scheme of $C$ over $K$, with connected component $\Aut^0_{C/K} \cong E$ and component group scheme $Q_C := \pi^0(\Aut_{C/K}) \cong \Aut_{C/K}/E$.
\item $p(C)$ is the \emph{period} of the $E$-torsor $C$, i.e., the order of the corresponding class in the Weil--Ch\^atelet group $H^1_{\et}(\Spec(K),E)$ of $E$.
\item $i(C)$ is the \emph{index} of $C$, i.e., the greatest common divisor of all degrees of field extensions of $K$ over which $C$ acquires a rational point.
\item $\tau$ (resp. $\tau_d$) denotes the natural action of $\Aut^0_{C/K} \cong E$ on $\Pic_{C/K}$ (resp. $\Pic^d_{C/K}$).
\end{itemize}

%
%


By definition, $\tau_1$ is just the restriction of the group law $m$ on $\Pic_{C/K}$ to $E \times C$. The other $\tau_d$ can be described similarly as follows.
\begin{Lemma} \label{compareaction}
For $d \in \mathbb{Z}$, let $[d]: \Pic_{C/K} \to \Pic_{C/K}$ be multiplication by $d$. Then, $\tau_d = m\circ ([d] \times {\rm id})$ holds. In particular, $\tau_0$ is the trivial action.
\end{Lemma}

\prf
The action $\tau$ is compatible with multiplication $m$ and inversion $i$ on $\Pic_{C/K}$. Hence, the following two diagrams commute for all $n,n' \in \ZZ$:

\makebox[\textwidth]{%
\xymatrix{
\Pic_{C/K}^0 \times \Pic_{C/K}^n \times \Pic_{C/K}^{n'} \ar[d]^-{\tau_n \times \tau_{n'}} \ar[r]^-{{\rm id} \times m} & \Pic_{C/K}^0 \times  \Pic_{C/K}^{n+n'} \ar[d]^-{\tau_{n+n'}}  \\
\Pic_{C/K}^n \times \Pic_{C/K}^{n'} \ar[r]^-{m} & \Pic_{C/K}^{n+n'} 
}
}

\bigskip and

\makebox[\textwidth]{%
\xymatrix{
\Pic_{C/K}^0 \times \Pic_{C/K}^n \ar[d]^-{\tau_n} \ar[r]^-{{\rm id} \times i} & \Pic_{C/K}^0 \times  \Pic_{C/K}^{-n} \ar[d]^-{\tau_{-n}}  \\
\Pic_{C/K}^n \ar[r]^-{i} & \Pic_{C/K}^{-n} 
}
}
The statement is true for $[d] = 1$.
Let $n,n' \in \ZZ$ and assume the statement holds for $d = n$ and $d = n'$. 
Then,
\begin{eqnarray*}
\tau_{n+n'} \circ ({\rm id} \times m) &=& m \circ (\tau_n \times \tau_{n'})  \\
&=& m \circ ((m \circ ([n] \times {\rm id})) \times (m \circ ([n'] \times {\rm id}))) \\
&=& m \circ ([n + n'] \times {\rm id}) \circ ({\rm id} \times m)
\end{eqnarray*}
and similarly
\begin{eqnarray*}
\tau_{-n} \circ ({\rm id} \times i) &=& m \circ ([-n] \times {\rm id}) \circ ({\rm id} \times i).
\end{eqnarray*}
Since ${\rm id} \times m$ and ${\rm id} \times i$ are epimorphisms, the claim follows for $d = n + n'$ and $d = -n$. Hence, the lemma is true for arbitrary $d \in \ZZ$.
\qed

\bigskip

\begin{Corollary} \label{periodaction}
The group $E(K)$ acts trivially on the quotient group $\Pic_{C/K}(K)/(p(C) \cdot \Pic_{C/K}(K))$.
\end{Corollary}

\prf
By \cite[Lemma 1]{Lichtenbaum}, the set $\Pic^d_{C/K}(K)$ is non-empty if and only if $p(C)$ divides $d$. On the other hand, for every $L \in E(K) = \Pic^0_{C/K}(K)$ and $M \in \Pic^d_{C/K}(K)$ we have $\tau(L,M) = \tau_d(L,M) = M + dL$ by the previous lemma. Hence, $E(K)$ acts trivially on $\Pic_{C/K}(K)/(p(C) \cdot \Pic_{C/K}(K))$.
\qed
\bigskip

Since $E \cong \Aut^0_{C/K}$ acts trivially (via $\tau_0$) on $E \cong \Pic^0_{C/K}$, there is an induced action of $Q_C =  \Aut_{C/K}/E$ on $E$.
%
The following lemma shows that the group of rational points $Q_C(K)$ tends to be rather small.

\begin{Lemma}\label{QC}
The action of $Q_C$ on $E$ is faithful and if $g \in Q_C(K)$, then  $\ord(g) \in \{1,2,3,4,6\}$. Moreover,
\begin{enumerate}[(i)]
\item if $\ord(g) = 6$, then $i(C) = 1$.
\item if $\ord(g) = 4$, then $i(C) \mid 2$.
\item if $\ord(g) = 3$, then $i(C) \mid 3$.
\item if $\ord(g) = 2$, then $i(C) \mid 4$.
\end{enumerate}
If $\ord(g)$ is prime to $\Char(K)$, then
\begin{enumerate}[(a)]
\item if $\ord(g) = 4$, then $\ZZ/2\ZZ \subseteq E(K)$.
\item if $\ord(g) = 3$, then $E$ contains a twisted form of the constant group scheme $\underline{\ZZ/3\ZZ}$.
\end{enumerate}
\end{Lemma}

\prf
First, assume that $K$ is algebraically closed. In this case, there is a point $P \in C(K)$ and any such point provides us with a splitting of the connected-\'etale sequence for $\Aut_{C/K}$ by identifying $Q_C$ with the subgroup scheme $\Aut_{C/K, P} \subseteq \Aut_{C/K}$ of automorphisms fixing $P$. Moreover, translation by $P$ induces a $Q_C$-equivariant isomorphism $E = \Pic^0_{C/K} \to \Pic^1_{C/K} = C$. In particular, the action of $Q_C$ on $E$ is faithful. Now, it is well known (see e.g. \cite{Silverman}) that if $g$ is non-trivial, then $\ord(g) = 2,3,4$ resp. $6$ and the fixed locus ${\rm Fix}(g) \subseteq E$ of $g$ is a closed subscheme of $E$ of length $4,3,2$ resp. $1$. Moreover, if $\ord(g)$ is prime to $\Char(K)$, then ${\rm Fix}(g)(K) \subseteq E(K)$ is a subgroup isomorphic to $\ZZ/2\ZZ \times \ZZ/2\ZZ$, $\ZZ/3\ZZ$, $\ZZ/2\ZZ$ resp. $\{1 \}$.

This shows that, over an arbitrary field, the action of $Q_C$ on $E$ is faithful and $\ord(g) \in \{1,2,3,4,6\}$. Moreover,
${\rm Fix}(g) \subseteq E$ is a closed subgroup scheme of length $1,2,3$ resp. $4$ on $C$ if $g$ is of order $6,4,3$ resp. $2$. Since $g$ is defined over $K$, so is ${\rm Fix}(g)$ and we deduce that $i(C)$ divides the degree of ${\rm Fix}(g)$, because $C$ acquires a rational point when base changed to ${\rm Fix}(g)$. This yields claims (i) to (iv).

If $\ord(g)$ is prime to $\Char(K)$, then ${\rm Fix}(g) \subseteq E$ is an \'etale subgroup scheme which becomes isomorphic to the constant group scheme $\underline{\ZZ/2\ZZ}$ resp. $\underline{\ZZ/3\ZZ}$ if $\ord(g) =4$ resp. $\ord(g) = 3$ over $\overline{K}$, i.e., ${\rm Fix}(g)$ is a twisted form of these groups. Since $\ZZ/2\ZZ$ admits no non-trivial automorphisms, there is no non-trivial twisted form of $\underline{\ZZ/2\ZZ}$, hence $\ZZ/2\ZZ \subseteq E(K)$.
\qed
\bigskip


After having studied automorphisms of $C$ as a $K$-scheme, we can also study equivariant automorphisms in the following sense.

\begin{Definition}\label{equivariant}
For a morphism of schemes $f: S \to T$, we define the automorphism group $\Aut(f) \subseteq \Aut(S) \times \Aut(T)$ of the morphism $f$ to be the subgroup of pairs $(g,h)$ such that the following diagram commutes

\makebox[\textwidth]{%
\xymatrix{
S \ar[r]^-{g} \ar[d]_f & S \ar[d]^f \\
T  \ar[r]^-{h} & T
}
}

\noindent The group $\Aut(f)$ comes with two projections which we will denote by $p_S: \Aut(f) \to \Aut(S)$ and $p_T: \Aut(f) \to \Aut(T)$.
\end{Definition}

We have the following lemma for a smooth curve $c:C \to \Spec(K)$ of genus one  over a field $K$ with Jacobian $e:E \to \Spec(K)$, which shows how equivariant automorphisms of $C$ act on $E$.

\begin{Lemma}\label{inducedactionpicardgenerically}
There is a morphism $\varphi: \Aut(c) \to \Aut(e)$ such that
\begin{enumerate}[(i)]
\item $p_K \circ \varphi = p_K$,
\item ${\rm Im}(\varphi)$ preserves $\mathcal{O}_C$, and
\item $E(K) = {\rm Ker}(\varphi) \subseteq {\rm Ker}(p_K)$.
\end{enumerate}
\end{Lemma}

\prf
Existence of a morphism $\varphi$ with properties $(i)$ and $(ii)$ follows from the fact that the construction of the Picard scheme over $K$ is functorial and commutes with changing the base (see \cite[Exercise 4.4]{Kleiman}). By $(i)$, we have ${\rm Ker}(\varphi) \subseteq {\rm Ker}(p_K) = \Aut(C) = \Aut_{C/K}(K)$, and $E(K) \subseteq {\rm Ker}(\varphi)$ follows from Lemma \ref{compareaction}. On the other hand, by Lemma \ref{QC}, $Q_C(K)$ acts faithfully on $E$, hence ${\rm Ker}(\varphi) \subseteq E(K)$.
\qed
\subsection{Translation into the setting of elliptic surfaces over $k$}
We fix the following setup.

\begin{itemize}
\item $f:X \to B$ is an elliptic surface over an algebraically closed field $k$.
\item $K = k(B)$ is the function field of $B$.
\item $C$ is the generic fiber of $f$.
\item $E$ is the Jacobian of $C$.
\item $J(f): J(X) \to B$ is the Jacobian fibration of $f$, which is the minimal proper smooth model of its generic fiber, the Jacobian $E$ of $C$. We denote its zero section by $\sigma$.
\item ${\rm MW}(J(f)) := E(K)$ is the Mordell--Weil group of $J(f)$ which comes with a natural action on $X$
$$
\tau_1(K): {\rm MW}(J(f)) \times X \to X.
$$
\end{itemize}

First, note that since the Brauer group $\Br(K)$ of $K$ is trivial by Tsen's Theorem \cite{Tsen}, period and index coincide in this setting by a result of S. Lichtenbaum.

\begin{Theorem}[\cite{Lichtenbaum}]\label{periodindex}
In the above setting, we have $p(C) = i(C)$.
\end{Theorem}

Moreover, there is a split surjective restriction morphism
\begin{equation*}
\rho: \Pic(X) \to \Pic(C).
\end{equation*}
The section is induced by the map sending an effective divisor on $C$ to its closure in $X$. The kernel of $\rho$ consists of divisors whose support is contained in finitely many fibers of $f$. Denote the kernel of $\rho$ by $\Pic^{\rm v}(X)$ and the image of the section of $\rho$ by $\Pic^{\rm h}(X)$. 
%
Our study of automorphisms of curves of genus one translates to this situation as follows.

\begin{Corollary}\label{actionhorizontalpicard}
The action on $\Pic^{\rm h}(X)/(i(C)\Pic^{\rm h}(X))$ induced by $\tau_1(K)$ is trivial. If, moreover, all fibers of $f$ are irreducible, then $\tau_1(K)$ induces the trivial action on $\Pic(X)/(i(C)\Pic^{\rm h}(X))$.
\end{Corollary}

\prf
We have natural isomorphisms $\Pic^{\rm h}(X) \cong \Pic(C) \cong \Pic_{C/K}(K)$, where the second isomorphism follows from $\Br(K)= 0$ using the Hochschild-Serre spectral sequence (see \cite{Kleiman}). Now, the first claim follows immediately from Corollary \ref{periodaction} and Theorem \ref{periodindex}.

Suppose that all fibers of $f$ are irreducible. Then, $\Pic^{\rm v}(X) \cong \ZZ$ is generated by a rational multiple of the class of a fiber of $f$. Since $\tau_1(K)$ preserves $f$, the induced action on $\Pic^{\rm v}(X)$ is trivial. Moreover, $\tau_1(K)$ maps $\Pic^{\rm h}(X)$ to $\Pic^{\rm h}(X)$, hence it acts trivially on $\Pic(X)/(i(C)\Pic^{\rm h}(X))$ by the first claim.
\qed

\begin{Corollary}\label{inducedactionpicard}
With notation as in Definition \ref{equivariant}, there is a morphism $\varphi: \Aut(f) \to \Aut(J(f))$ such that
\begin{enumerate}[(i)]
\item $p_B \circ \varphi = p_B$,
\item ${\rm Im}(\varphi)$ preserves $\sigma$, and
\item ${\rm MW}(J(f)) = {\rm Ker}(\varphi)$.
\end{enumerate}
\end{Corollary}

\prf
This follows immediately from Lemma \ref{inducedactionpicardgenerically}, since every automorphism of $C \to \Spec(K)$ (resp. $E \to \Spec(K)$) extends uniquely to an automorphism of $f$ (resp. $J(f)$), because both of them are the unique relatively minimal proper smooth model of their generic fiber.
\qed

\section{Rational elliptic surfaces with two marked fibers} \label{section3}
Throughout this section, we assume $\Char(k) \neq 2$. Here, we collect results on rational elliptic surfaces which appear as Jacobians of elliptic fibrations on unnodal Enriques surfaces. In particular, the base of these rational elliptic fibrations comes with two marked points corresponding to the two double fibers of the elliptic fibration of the Enriques surface.

We consider triples $(f,\sigma,\{p_1,p_2\})$, called \emph{2-marked rational elliptic surfaces}, where $f: J \to \PP^1$ is a rational elliptic surface, $\sigma: \PP^1 \to J$ is a section of $f$ and $p_1,p_2 \in \PP^1$ are closed points. An isomorphism between $(f,\sigma,\{p_1,p_2\})$ and $(f',\sigma',\{p_1',p_2'\})$ is an isomorphism between $f$ and $f'$ compatible with sections and the sets of marked points. 

Let ${\rm Art}_k$ be the category of Artinian local $k$-algebras with residue field $k$. A deformation of a $2$-marked rational elliptic surface $(f,\sigma,\{p_1,p_2\})$ over $A \in {\rm Art}_k$ is a commutative diagram of $A$-flat schemes
$$
\xymatrix{ \mathcal{J} \ar[r]^{\overline{f}} \ar[dr] & \mathcal{C} \ar[d] \ar@/_1.5pc/[l]_{\overline{\sigma}} \\
& {\rm Spec}(A) \ar@<-.5ex>@/_1.5pc/[u]_{\overline{p}_2} \ar@<.5ex>@/_1.5pc/[u]^{\overline{p}_1}
}
$$
together with an isomorphism identifying the restriction of the above diagram to the special fiber with $(f,\sigma,\{p_1,p_2\})$. Similary, we can define families of such triples and isomorphisms of deformations and families.
We say that a family of $2$-marked rational elliptic surfaces $(\overline{f},\overline{\sigma},\{\overline{p}_1,\overline{p}_2\})$ over an algebraic space $\mathcal{S}$ is universal at $s \in \mathcal{S}(k)$, if its restriction to $\hat{\mathcal{O}}_{S,s}$ induces the universal formal deformation of the fiber $(\overline{f}_s,\overline{\sigma}_s,\{(\overline{p}_1)_s,(\overline{p}_2)_s\})$.

Via the natural action of ${\rm PGL}_2$ on $\PP^1$, every $2$-marked rational surface $(f,\sigma,\{p_1,p_2\})$ is isomorphic to $(f,\sigma,\{0,\infty\})$, and similarly for its deformations, where $0$ and $\infty$ are defined with respect to some fixed coordinate $t$ on $\PP^1$. Hence, we can sometimes restrict to triples of the form $(f,\sigma,\{0,\infty\})$.

In our case, these triples will arise as Jacobians of elliptic fibrations on unnodal Enriques surfaces and the two marked points will correspond to the double fibers of the elliptic fibration of the Enriques surface. Therefore, by the facts recalled in Section \ref{EllfibronEnriques}, they will satisfy the following extra conditions.

\begin{Definition}
We say that a $2$-marked rational elliptic surface $(f,\sigma,\{p_1,p_2\})$ satisfies $(\ast)$ if
\begin{itemize}
\item all fibers of $f$ are irreducible, and
\item the fibers of $f$ over $p_1$ and $p_2$ are of type $\I_0$ or $\I_1$. In particular, the discriminant $\Delta$ of $f$ vanishes at most once at $p_1$ and $p_2$.
\end{itemize}
\end{Definition}

\begin{WeierstrassModelx} Note that for a family of rational elliptic surfaces $(f: J \to \PP^1_A,\sigma,\{p_1,p_2\})$ over $A \in {\rm Art}_k$ , we have $R^1f_*\mathcal{O}_J = \mathcal{O}_{\PP^1_A}(-1)$ by the canonical bundle formula. By the classical theory of Weierstra{\ss} models (see e.g. \cite{MumfordSuominen}, \cite{Lang3} for the case of characteristic $3$, and also \cite{Seiler}), a family of rational elliptic surfaces with only irreducible fibers over $A \in {\rm Art}_k$ is isomorphic to its Weierstra{\ss} model, which is a closed subscheme of the projective bundle $\PP(\mathcal{O}_{\PP^1_A} \oplus \mathcal{O}_{\PP^1_A}(-2) \oplus \mathcal{O}_{\PP^1_A}(-3))$ over $\PP^1_A$ given by an affine equation of the form
\begin{equation} \label{Weierstrassequation}
y^2 = x^3 + a_2x^2 + a_4x +a_6,
\end{equation}
with $a_i = \sum_{j = 0}^i a_{ij}t^j$ and $a_{ij} \in A$. Such an affine Weierstra{\ss} equation in fact determines a family of $2$-marked rational elliptic surfaces by choosing the fibers over $0$ and $\infty$ as the marked fibers and the section $\sigma$ as the section given by $z = 0$ in the homogenization of Equation (\ref{Weierstrassequation}).
\end{WeierstrassModelx}

\begin{Remark}
The above situation is much simpler than the analogous one for deformations of arbitrary elliptic surfaces (e.g. as in \cite{Seiler}), since $\PP^1$ and all line bundles on it are infinitesimally rigid, so only the $a_i$ can deform non-trivially. Also, the vector bundle $V$ occuring in \cite[Theorem 1.1 (B) (2)]{Lang3} is isomorphic to $\mathcal{O}_{\PP^1_A}(-2) \oplus \mathcal{O}_{\PP^1_A}$, since
$
{\rm Ext}^1_{\mathcal{O}_{\PP^1}}(\mathcal{O}_{\PP^1}(-2),\mathcal{O}_{\PP^1}) = 0.
$
\end{Remark}

The $a_i$ are not uniquely determined by $(f,\sigma,\{p_1,p_2\})$, but the ambiguity is easy to control.

\begin{Lemma} \label{substitutions}
Let $(f,\sigma,\{0,\infty\})$ resp. $(f',\sigma',\{0,\infty \})$ be families of $2$-marked rational elliptic surfaces over $A \in {\rm Art}_k$ given by equations
$y^2 = x^3 + a_{2}x^2 + a_4x + a_6$ resp. $(y')^2 = (x')^3 + a_2'(x')^2 + a_4'(x') + a_6'$. Then, any isomorphism $g$ from $(f,\sigma,\{0,\infty \})$ to $(f',\sigma', \{0,\infty \})$ is given by one of the following substitutions
\begin{enumerate}[(i)]
\item 
\begin{eqnarray*}
x' &\mapsto& u^{-2} x + r \\
y' &\mapsto& u^{-3} y  \\
t  &\mapsto& \lambda t
\end{eqnarray*}
with $\lambda, u \in A^\times$ and $r \in A[t]$ such that
\begin{eqnarray*}
a_2 &=& u^2(a_2'(\lambda t) + 3r) \\
a_4 &=& u^4(a_4'(\lambda t) + 2a_2'(\lambda t)r + 3r^2) \\
a_6 &=& u^6(a_6'(\lambda t) + a_4'(\lambda t)r + a_2'(\lambda t)r^2 + r^3),
\end{eqnarray*}or
\item
\begin{eqnarray*}
x' &\mapsto& (u t)^{-2} x + r \\
y' &\mapsto& (u t)^{-3} y  \\
t  &\mapsto& (\lambda t)^{-1}
\end{eqnarray*}
with $\lambda, u \in A^\times$ and $r \in A[t^{-1}]$ such that
\begin{eqnarray*}
a_2 &=& u^2t^2 \left(a_2'\left( (\lambda t)^{-1}\right) + 3r \right) \\
a_4 &=& u^4t^4 \left(a_4'\left( (\lambda t)^{-1}\right) + 2a_2' \left( (\lambda t)^{-1} \right)r + 3r^2 \right) \\
a_6 &=& u^6t^6 \left(a_6'\left( (\lambda t)^{-1} \right) + a_4'\left( (\lambda t)^{-1} \right)r + a_2'\left( (\lambda t)^{-1} \right)r^2 + r^3 \right).
\end{eqnarray*}
\end{enumerate}
We write $g = (h,u,r)$, where $h$ is the action of $g$ on $\PP^1$.
\end{Lemma}

\prf
If $g$ fixes $\PP^1$, then this is a standard exercise and explained for example in \cite{Tate}. If $g$ does not fix $\PP^1$, then $g$ acts on $t$ via an automorphism $h$ of $\PP^1_A$ preserving $\{0,\infty\}$. Thus, $h$ is of the form $t \mapsto \lambda t$ or $t \mapsto (\lambda t)^{-1}$. Since $g$ induces an isomorphism between $(f,\sigma,\{0,\infty \})$ and the pullback of $(f',\sigma',\{0,\infty \})$ along $h$, we reduce to the case where $g$ fixes $\PP^1$.
\qed

\smallskip
Using this explicit description of isomorphisms between $2$-marked rational elliptic surfaces, we can write down normal forms and universal deformations for them.

\begin{Lemma} \label{Weierstrassuni}
Let $(f,\sigma,\{0,\infty\} \})$ be a $2$-marked rational elliptic surface satisfying $(\ast)$. Consider the family of $2$-marked rational elliptic surfaces $(g,\tau,\{0,\infty\})$ defined by $y^2 = x^3 + a_2x^2 + a_4x + a_6$ over $\mathbb{A}^{15}  = \Spec(k[\{a_{ij}\}_{i \in \{2,4,6\},0 \leq j \leq i}])$

\begin{enumerate}[(i)]
\item If $\Char(k) \neq 3$, then the restriction of the family $(g,\tau,\{0,\infty\})$ to at least one of the following three closed subsets of $\mathbb{A}^{15}$ induces the universal formal deformation of $(f,\sigma,\{0,\infty\} \})$
\begin{itemize}
\item $V(a_{20},a_{21},a_{22},a_{40} - 1,a_{60} - 1)$
\item $V(a_{20},a_{21},a_{22},a_{40} - 1,a_{66} - 1)$
\item $V(a_{20},a_{21},a_{22},a_{60} - 1,a_{66} - 1)$
\end{itemize}
\item If $\Char(k) = 3$, then the restriction of the family $(g,\tau,\{0,\infty\})$ to at least one of the following three closed subsets of $\mathbb{A}^{15}$ induces the universal formal deformation of $(f,\sigma,\{0,\infty\} \})$
\begin{itemize}
\item $V(a_{40},a_{41},a_{42},a_{20} - 1,a_{22} - 1)$
\item $V(a_{40},a_{41},a_{42},a_{20} - 1,a_{44} - 1)$
\item $V(a_{60},a_{61},a_{62},a_{40} - 1,a_{44} - 1)$
\end{itemize}
\end{enumerate}
\end{Lemma}

\prf
We have to prove that every deformation $(f',\sigma',\{0,\infty\})$ of $(f,\sigma,\{0,\infty\})$ over $A \in {\rm Art}_k$ arises uniquely as a pullback of one of the given families. Denote the maximal ideal of $A$ by $\mathfrak{m}$ and assume that $(f',\sigma',\{0,\infty\})$ is given by an affine Weierstra{\ss} equation as in Equation (\ref{Weierstrassequation}).

First, assume $\Char(k) \neq 3$. Then, we can apply the substitution $x \mapsto x - \frac{1}{3}a_2$ to achieve $a_2 = 0$. The discriminant of $f'$ is then given by $\Delta = -16(4a_4^3 + 27a_6^2)$. Since $(f,\sigma,\{0,\infty\})$ satisfies $(\ast)$, either $a_{40}$ or $a_{60}$ and either $a_{44}$ or $a_{66}$ are units. Thus, after rescaling $x$ and $y$ as in Lemma \ref{substitutions} and possibly interchanging $t$ with $1/t$, we can assume $(a_{40},a_{60}) = 1$ or $(a_{40},a_{66}) = 1$ or $(a_{60},a_{66}) = 1$. This shows that $(f',\sigma',\{0,\infty\})$ arises as pullback of one of these three families. Now, for the uniqueness statement, it remains to show that any isomorphism of equations of this form is non-trivial modulo $\mathfrak{m}$. But by Lemma \ref{substitutions} and with the same notation, a substitution that is trivial mod $\mathfrak{m}$ acts on $\PP^1$ as $h: t \mapsto \lambda t$. Moreover, it satisfies $\lambda = u = 1~{\rm mod}~\mathfrak{m}$. On the other hand, to preserve the form of the equation, we must have $\lambda^{12} = u^{12} = 1$. Since $\Char(k) \neq 2,3$, this implies $\lambda = u = 1$ and proves the claim.

Now, assume $\Char(k) = 3$. The discriminant of $f'$ is given by $\Delta = a_2^2a_4^2 - a_2^3a_6 - a_4^3$. Since $(f,\sigma,\{0,\infty\})$ satisfies $(\ast)$, either $a_{20}$ or $a_{40}$ and either $a_{22}$ or $a_{44}$ are units. Similarly to the case $\Char(k) \neq 3$, we can thus assume $(a_{20},a_{22}) = 1$ or $(a_{20},a_{44}) = 1$ or $(a_{40},a_{44}) = 1$.
In the first two cases, we apply a substitution of the form $x \mapsto x + r$ to set $a_{40} = a_{41} = a_{42} = 0$. An isomorphism of equations satisfying the five conditions that is trivial mod $\mathfrak{m}$ has to satisfy $\lambda = u = 1~{\rm mod}~\mathfrak{m}$ and $r = 0~{\rm mod}~\mathfrak{m}$. On the other hand, we need $u^2 = 1$ since $a_{20} = 1$, $\lambda^4 = 1$ since $a_{22} = 1$ or $a_{44} = 1$ and $r = 0$ since $a_{40} = a_{41} = a_{42}$ and $r$ is of degree at most $2$. Thus, $(u,\lambda,r) = (1,1,0)$. In the third case, we can apply a substitution of the form $x \mapsto x + r$ to set $a_{60} = a_{61} = a_{62} = 0$ and the uniqueness is proved using the same method as in the other cases.
\qed
\smallskip

Next, we will study automorphisms of $2$-marked rational elliptic surfaces $(f,\sigma,\{0,\infty\})$ satisfying $(\ast)$. We will see in the next section, that extra automorphisms of Enriques surfaces induce automorphisms of order $2^n$ in $\Aut(f,\sigma,\{0,\infty\})$. Therefore, let us define the sets
\begin{eqnarray*}
\Aut_2(f,\sigma,\{0,\infty\}) &:=& \{g \in \Aut(f,\sigma,\{0,\infty\}) ~|~ \ord(g) = 2^n \text{ for some } n \in \ZZ \} \\
\overline{\Aut}_2(f,\sigma,\{0,\infty\}) &:=& {\rm Im}(\Aut_2(f,\sigma,\{0,\infty\}) \to \Aut(\PP^1))
\end{eqnarray*}

%


\noindent In the following lemma, we study the set $\Aut_2(f,\sigma,\{0,\infty\})$.

\begin{Lemma} \label{Weierstrassmain}
Let $(f,\sigma,\{0,\infty\})$ be a $2$-marked rational elliptic surface satisfying $(\ast)$ given by an affine Weierstra{\ss} equation as in Equation (\ref{Weierstrassequation}). If $\Char(k) \neq 3$, we further assume $a_2= 0$ and if $\Char(k) = 3$, we assume that either $a_{20} = 1, a_{40} = a_{41} = a_{42} = 0$ or $a_{40} = 1, a_{60} = a_{61} = a_{62} = 0$.

Let $g = (h,u,r) \in \Aut_2(f,\sigma,\{0,\infty\})$ with ${\rm id} \neq h$. Then, $h,u,r,a_2,a_4$ and $a_6$ are as in Tables \ref{Table1}, \ref{Table2} and \ref{Table3}, where $\lambda \in k^\times$, $r_0,r_1 \in k$ and $\zeta_4$ is a primitive $4$-th root of unity:

\begin{table}[h!]

\centering
\resizebox{\columnwidth}{!}{%
$
\begin{array}{|l|l|l|l|l|} \hline
h & u^2 & r& a_4 & a_6 \\ \hline
t \mapsto -t & 1 &0  & a_{40} + a_{42}t^2 + a_{44}t^4 & a_{60} + a_{62}t^2 + a_{64}t^4 + a_{66}t^6 \\ 
t \mapsto -t & -1 &0  & a_{40} + a_{42}t^2 + a_{44}t^4 & a_{61}t + a_{63}t^3 + a_{65}t^5 \\ 
t \mapsto \zeta_4 t& 1 &0  & a_{40} + a_{44}t^4 & a_{60} + a_{64} t^4 \\
t \mapsto \zeta_4 t & -1  &0  & a_{40} + a_{44}t^4 & a_{62}t^2 + a_{66} t^6 \\ \hline
t \mapsto (\lambda t)^{-1} & \lambda &0  & a_{40} + a_{41}t + a_{42}t^2 +  \lambda a_{41}t^3 + \lambda^2 a_{40}t^4 & a_{60} + a_{61}t +a_{62}t^2 + a_{63}t^3 + \lambda a_{62}t^4 + \lambda^2 a_{61}t^5 + \lambda^3  a_{60}t^6 \\
t \mapsto  (\lambda t)^{-1} & - \lambda &0   & a_{40} + a_{41}t + a_{42}t^2 +  \lambda a_{41}t^3 + \lambda^2 a_{40}t^4 & 
a_{60} + a_{61}t +a_{62}t^2 - \lambda a_{62}t^4 - \lambda^2 a_{61}t^5 - \lambda^3  a_{60}t^6 \\ \hline
\end{array}
$
}
\vspace{1mm}
\caption{$2$-marked rational elliptic surfaces with extra automorphisms in $\Char(k) \neq 3$}
\label{Table1}
\end{table}
\vspace{-6mm}
\begin{table}[h!]
\centering
\resizebox{\columnwidth}{!}{%
$
\begin{array}{|l|l|l|l|l|l|} \hline
h & u^2 & r & a_2 & a_4 & a_6 \\ \hline
t \mapsto -t & 1 & 0 & 1 + a_{22}t^2 & a_{44}t^4 & a_{60} + a_{62}t^2  + a_{64}t^4 + a_{66}t^6 \\ 
t \mapsto \zeta_4 t & 1 & 0  & 1 & a_{44}t^4 & a_{60} + a_{64}t^4 \\ \hline
t \mapsto (\lambda t)^{-1} & \lambda   & r_0 -  \frac{r_0}{\lambda} t^{-2} &  1 + a_{21}t + \lambda t^2 & -\lambda a_{21} r_0 t^3 - \lambda^2 r_0 t^4 & \begin{array}{l}
a_{60} + a_{61}t +a_{62}t^2 + a_{63}t^3   \\ + \lambda (a_{62} - \lambda r_0^2)t^4 \\ + \lambda^2(a_{61} + a_{21}r_0^2)t^5  \\ + \lambda^3 (a_{60} + r_0^2 + r_0^3)t^6
\end{array}
\\ \hline
t \mapsto (\lambda t)^{-1} & - \lambda   & r_0 + r_1 t^{-1} + \frac{r_0}{\lambda}t^{-2} &  1 - \lambda t^2 & -\lambda^2(r_1 t^3 +r_0 t^4) & 
\begin{array}{l}
a_{60} + a_{61}t +a_{62}t^2  + (\lambda^3 r_1^3 - \lambda^2 r_0 r_1) t^3  \\ - \lambda(a_{62} + \lambda r_0^2 + \lambda^2 r_1^2) t^4 \\ - \lambda^2(a_{61} - \lambda r_0r_1) t^5  \\ - \lambda^3(a_{60} + r_0^2 + r_0^3) t^6
\end{array}
\\


 \hline
\end{array}
$
}
\vspace{1mm}
\caption{$2$-marked rational elliptic surfaces with extra automorphisms in $\Char(k) = 3$ with $a_{20} = 1, a_{40} = a_{41} = a_{42} = 0$}
\label{Table2}
\end{table}
\vspace{-6mm}
\begin{table}[h!]
\centering
\resizebox{\columnwidth}{!}{%
$
\begin{array}{|l|l|l|l|l|l|} \hline
h & u^2 & r & a_2 & a_4 & a_6 \\ \hline
t \mapsto -t & 1 & 0 & a_{20} + a_{22}t^2 & 1 + a_{42}t^2 + a_{44}t^4 & a_{64}t^4 + a_{66}t^6 \\
t \mapsto -t & -1 & 0 \text{ or } \zeta_4 (1 - a_{42}t^2) & a_{21}t & 1 + a_{42}t^2 + a_{44}t^4 - ra_{21}t& a_{63}t^3 + a_{65}t^5 + (a_{44} - a_{42}^2)rt^4 \\
t \mapsto \zeta_4 t & 1 & 0 & a_{20} & 1 + a_{44}t^4 & a_{64}t^4 \\
t \mapsto  \zeta_4 t & -1 & 0 \text{ or } \pm \zeta_4 & a_{22}t^2 & 1 + a_{44}t^4 - ra_{22}t^2 & a_{66}t^6 + a_{44}rt^4 \\
 \hline
 
 t \mapsto (\lambda t)^{-1} & \lambda  & r_0 -  \frac{r_0}{\lambda} t^{-2} & 
 \begin{array}{l}
 a_{20} + a_{21}t \\ + \lambda a_{20}t^2 
 \end{array}
 &
   \begin{array}{l} 1 + a_{41}t + a_{42}t^2 \\ + \lambda(a_{41} - r_0 a_{21})t^3 \\ + \lambda^2(1 - r_0 a_{20})t^4
   \end{array} & 
\begin{array}{l}
a_{63}t^3 \\ + \lambda r_0 (a_{42} - \lambda - \lambda r_0 a_{20} )t^4 \\ +  \lambda^2 r_0 (a_{41} + r_0 a_{21})t^5 \\ + \lambda^3 r_0 (1 + a_{20}r_0 + r_0^2)t^6
\end{array}

     \\ \hline

 t \mapsto (\lambda t)^{-1} & - \lambda & r_0 + r_1 t^{-1} + \frac{r_0}{\lambda}t^{-2} & a_{20}  -\lambda a_{20}t^2 &
 \begin{array}{l}
   1 + a_{41}t + a_{42}t^2 \\ + \lambda(a_{41} - \lambda r_1 a_{20})t^3 \\ + \lambda^2(1 - r_0 a_{20})t^4
  \end{array}
    
    & \begin{array}{l} (\lambda^3 r_1^3 + \lambda a_{42}r_1 - \lambda a_{41}r_0 - \lambda^2 a_{20} r_0 r_1)t^3 \\ - \lambda (\lambda r_0 + \lambda a_{41}r_1 + a_{42}r_0 + \lambda a_{20}r_0^2 + \lambda^2 a_{20}r_1^2)t^4 \\ - \lambda^2(\lambda r_1 + a_{41}r_0 - \lambda a_{20}r_0r_1)t^5 \\ - \lambda^3(r_0 + a_{20}r_0^2 + r_0^3)t^6
\end{array}     \\ \hline

\end{array}
$
}
\vspace{1mm}
\caption{$2$-marked rational elliptic surfaces with extra automorphisms in $\Char(k) = 3$ with $a_{40} = 1, a_{60} = a_{61} = a_{62} = 0$}
\label{Table3}
\end{table}

\end{Lemma}

\prf
First, note that $g^{\ord(h)}$ fixes $\PP^1$ pointwise. Since $f$ has no reducible fibers, it has no $2$-torsion section (see \cite[p.111]{ShiodaSchuett}) and hence by Lemma \ref{QC} the order $\ord(g^{\ord(h)})$ is not divisible by $4$. On the other hand, by assumption, $g$ has order $2^n$ for some $n$. This implies that $\ord(g^{\ord(h)}) \leq 2$.

If $h$ is of the form $t \mapsto (\lambda t)^{-1}$, then $g^2$ maps $y$ to $u^{-6} \lambda^3 y$. Since $g^2$ has order $1$ or $2$, this implies $u^2 = \pm \lambda$. On the other hand, $g^2$ maps $x$ to $x + u^{-2}\lambda^2t^2r(t) + r((\lambda t)^{-1})$ and since it has order $1$ or $2$, this shows $u^{-2}\lambda^2t^2r(t) + r((\lambda t)^{-1}) = 0$. In particular, $r$ is a polynomial in $t^{-1}$ of degree at most $2$ and of the form given in the tables. Of course, if $\Char(k) \neq 3$, then $a_2 = 0$ forces $r = 0$. Calculating the $a_i$ in each case is tedious but straightforward.

Now, assume that $h$ is of the form $t \mapsto \lambda t$. Let us first treat the case $\Char(k) \neq 3$ where we know that $r = 0$. There, if $a_{i0} \neq 0$, then $u^i = 1$ and if $a_{ii} \neq 0$, then $\lambda^i u^i = 1$. Going through the possible combinations, we easily see that $\lambda^4 = 1$ and $u^{12} = 1$. Since $g^4$ maps $x$ to $u^{8}x$ and $g^4$ has order at most $2$, we see that $u^{16} = 1$. Combining this with $u^{12} = 1$, we see that $u^4 = 1$, i.e. $u^2 = \pm 1$. Again, $a_2 = 0$ forces $r = 0$ and it is easy to deduce the table.

If $\Char(k) = 3$ and $a_{20} = 1$, then $u^2 = 1$. Since either $a_{22} \neq 0$ or $a_{44} \neq 0$, we have $\lambda^4 = 1$. Moreover, $a_{40} = a_{41} = a_{42} = 0$ forces $r = 0$. Again, it is straightforward to calculate the $a_i$ in this case.

If $\Char(k) = 3$ and $a_{40} = 1$, then $u^4 = 1$. As in the previous paragraph, this implies $\lambda^4 = 1$. If $u^2 = 1$, then $r = 0$, since $g$ maps $x$ to $x + r$ which has odd order if $r \neq 0$, contradicting our assumptions on $g$. If $u^2 = -1$, however, $r$ is not necessarily trivial. But in this case, we must have $a_{20} = 0$ and since $a_{60} = 0$ the equality $a_{60} = u^6(a_{60}(\lambda t) + a_{40}(\lambda t)r + a_{20}(\lambda t)r^2 + r^3)$ shows that the constant term $r_0$ of $r$ satisfies $r_0 + r_0^3 = 1$. This is the reason for the three possibilities for $r$ in the table in these cases. We leave the computation of the $a_i$ in each case to the reader.
\qed

\bigskip

The above calculations allow us to give a bound on the size of $2$-groups contained in $\overline{\Aut}_2(f,\sigma,\{0,\infty\})$. This is the bound that eventually leads to the corresponding bound in Theorem \ref{mainA}.

\begin{Corollary}\label{extrabound}
Let $(f,\sigma,\{0,\infty \})$ be a $2$-marked rational elliptic surface satisfying $(\ast)$ and let $\overline{G}_2$ be a $2$-group contained in $\overline{\Aut}_2(f,\sigma,\{0,\infty\})$. Then, we have $|\overline{G}_2| \in \{1,2,4\}$.
\end{Corollary}

\prf
By Lemma \ref{Weierstrassmain}, $\overline{G}_2$ contains no element of order $8$. Since $\overline{G}_2$ preserves $\{0,\infty\}$, we have $\overline{G}_2 \subseteq \GG_m \rtimes \ZZ/2\ZZ$ and thus we only have to show that $\overline{G}_2$ does not contain the dihedral group $D_8$. From the tables in Lemma \ref{Weierstrassmain}, we see that this would imply that $(f,\sigma,\{0,\infty \})$ is given by one of the following equations
\begin{itemize}
\item $y^2 = x^3 + (a_{40} + \lambda^2a_{40}t^4)x$ \quad if $\Char(k) \neq 3$.
\item $y^2 = x^3 + (1 + \lambda^2t^4)x - \lambda^2 r_0 t^4$ \quad with $r_0^3 + r_0 = 0$ if $\Char(k) = 3$.
\end{itemize}
However, these equations define singular surfaces contradicting our assumption that $(f,\sigma,\{0,\infty \})$ satisfies $(\ast)$.
\qed

\bigskip
Finally, the explicit calculations of this section allow us to determine the dimension of the moduli space of $2$-marked rational elliptic surfaces $(f,\sigma,\{0,\infty\})$ with non-trivial $\overline{\Aut}_2(f,\sigma,\{0,\infty\})$.

\begin{Theorem} \label{Weierstrassmoduli}
Let $(f,\sigma,\{p_1,p_2\})$ be a family of $2$-marked rational elliptic surfaces satisfying $(\ast)$ over an algebraic space $\mathcal{S}$ locally of finite type over ${\rm Spec}(k)$. 
We denote the geometric fiber of $(f,\sigma,\{p_1,p_2\})$ over a point $t \in |\mathcal{S}|$ by $(f_{\bar t},\sigma_{\bar{t}},\{(p_1)_{\bar t},(p_2)_{\bar t}\})$. 
For $j \in \{2,4\}$, define the set
$$
\mathcal{Z}_j = \{t \in \mathcal{S}(k) ~|~ \overline{\Aut}_2(f_{\bar t},\sigma_{\bar{t}},\{(p_1)_{\bar t},(p_2)_{\bar t}\}) \text{ contains a } 2\text{-group of order at least }  j \} \subseteq |\mathcal{S}|.
$$

Assume that $(f,\sigma,\{p_1,p_2\})$ is universal at a point $s \in \mathcal{S}(k)$. Then, there is an open neighborhood $s \in U \subseteq |\mathcal{S}|$ such that 
\begin{enumerate}[(i)]
\item $U$ has dimension $10$,
\item $\mathcal{Z}_2 \cap U \subseteq U$ is closed of dimension at most $5$, and
\item $\mathcal{Z}_4 \cap U \subseteq U$ is closed of dimension at most $2$.
\end{enumerate}
\end{Theorem}

\prf
Since $\mathcal{S}$ is locally of finite type over $k$ and universal at $s$, it is an algebraization of the universal formal deformation in the sense of M. Artin \cite{Artin}. By \cite[Theorem 1.7]{Artin}, such an algebraization is unique in an \'etale neighborhood of $s \in S$. Recall that \'etale morphisms are open of relative dimension $0$ and thus both claims in the theorem can be checked on such an \'etale neighborhood. Hence, it suffices to prove the theorem for one specific algebraization of the universal formal deformation of $(f_s,\sigma_s,\{(p_1)_s,(p_2)_s\})$.

In Lemma \ref{Weierstrassmain} we have described such algebraizations, all of which are $10$-dimensional. Using the explicit equations in Lemma \ref{Weierstrassmain}, we can calculate the sets $\mathcal{Z}_2$ and $\mathcal{Z}_4$ explicitly for the families described in Lemma \ref{Weierstrassuni}. There, $\mathcal{Z}_2$ and $\mathcal{Z}_4$ are closed and they have the claimed dimension. Since one of the three families in Lemma \ref{Weierstrassuni} induces the universal formal deformation of $(f_s,\sigma_s,\{(p_1)_s,(p_2)_s\})$, this proves the claim.
\qed

\bigskip

\section{Automorphisms of unnodal Enriques surfaces} \label{section4}
%
%

In order to be able to apply our knowledge of rational elliptic surfaces to an unnodal Enriques surface $X$ via the machinery of Section \ref{section2}, we have to prove that all extra automorphisms of $X$ are induced by automorphisms of $X$ preserving an elliptic fibration. The key step in this direction is the following theorem.

\begin{Theorem}\label{2group} Let $X$ be an unnodal Enriques surface.
Then, the extra automorphism group $\overline{\Aut}(X)$ is a $2$-group. 
\end{Theorem}

\prf
By Corollary \ref{possibleorders}, we know that the order of $\overline{\Aut}(X) \subseteq {\rm O}^+(10,\mathbb{F}_2)$ divides $2^{21} \cdot 3^5 \cdot 5^2 \cdot 7 \cdot 17 \cdot 31$, so it suffices to exclude that $\overline{\Aut}(X)$ contains elements of order $3,5,7,17$ or $31$. We let $\pi: Y \to X$ be the K3 cover of $X$.

\vspace{3mm}
\noindent \textbf{Claim 1:} If $k = \mathbb{C}$, then there is no element of order $17$ or $31$ in $\overline{\Aut}(X)$.

\indent \indent Using the notation of Lemma \ref{extendability}, we embed $S := \pi^*(\Num(X)) \subseteq \Num(Y) \hookrightarrow H^2(Y,\mathbb{Z}) =:  M$ and let $K = S^\perp$ and $T = \Num(Y)^\perp$. Then, $K = L \oplus T$ for some primitive and negative definite sublattice $L \subset H^2(Y,\mathbb{Z})$. 

Now, let $g$ be an automorphism of $X$ whose image $\overline{g} \in \overline{\Aut}(X)$ has order $17$ resp. $31$ and let $\tilde{g}$ be a lift of $g$ to $Y$. From $\tilde{g}^*(S) = S$ and  $\tilde{g}^*(\Num(Y)) = \Num(Y)$, we deduce that $\tilde{g}^*(T) = T$ and $\tilde{g}^*(L) = L$. Note that $\tilde{g}^*|_L$ has finite order, since $L$ is negative definite, and $\tilde{g}^*|_T$ has finite order, since $\tilde{g}^*$ preserves the period of $Y$. Therefore, $\tilde{g}^*|_K$ also has finite order and ${\rm rk}(K) = 22 - {\rm rk}(S) = 12$ implies that $\varphi(\ord(\tilde{g}^*|_K)) \leq 12$, where $\varphi$ is the Euler $\varphi$-function. On the other hand, by Lemma \ref{identification}, $\tilde{g}^*$ acts faithfully on the discriminant of $S$ and hence also on the discriminant of $K$. Therefore, $\ord(\tilde{g}^*|_K)$ is divisible by $17$ resp. $31$. This is a contradiction and proves Claim $1$.

\vspace{3mm}
\noindent \textbf{Claim 2:} If $\Char(k) > 2$ and $Y$ has finite height, then there is no element of order $17$ or $31$ in $\overline{\Aut}(X)$.

\indent \indent Again, let $g \in \Aut(X)$ be an automorphism whose image $\overline{g} \in \overline{\Aut}(X)$ has order $17$ resp. $31$ and let $\tilde{g}$ be a lift of $g$ to $Y$.
This time, we set $T_{crys} = (\Num(Y) \otimes W)^\perp \subseteq H^2_{crys}(Y/W)$. By \cite[Proposition 3.6]{Jang2}, the image of the representation $\chi_{crys,Y}: \Aut(Y) \to {\rm O}(T_{crys})$ is finite. Since ${\rm rk}(T_{crys}) = 22 - {\rm rk}(\Num(Y)) \leq 12$, we can assume $\chi_{crys,Y}(\tilde{g}) = {\rm id}$ without loss of generality. Hence, by \cite[Theorem 3.3]{Jang1} and its proof, there is a N\'eron--Severi group preserving lift $Y_\eta$ of $Y$ to characteristic $0$ that also lifts $\tilde{g}$ to an automorphism $\tilde{g}_\eta$ of $Y_\eta$. Since specialization of line bundles is equivariant with respect to the actions of the automorphism groups, the actions of $\tilde{g}_\eta$ and $\tilde{g}$ on $\Num(Y_\eta) \overset{\sim}{\to} \Num(Y)$ coincide. Now, we can exclude the existence of $\tilde{g}_\eta$ as in Claim $1$, thus proving the non-existence of $g$ and Claim $2$.

\vspace{3mm}
\noindent \textbf{Claim 3:} If $\Char(k) > 2$ and $Y$ is supersingular, then there is no element of order $17$ or $31$ in $\overline{\Aut}(X)$.

\indent \indent Let $g,\overline{g}$ and $\tilde{g}$ be as before. We set
$S = \pi^*(\Num(X)) \subset \Num(Y) =: M$ and $K := S^\perp$.
Since $Y$ is supersingular, the discriminant of $\Num(Y)$ is $-p^{2\sigma}$, where $p = \Char(k)$ and $\sigma$ is the Artin invariant of $Y$ (see \cite{Artin2}). With notation as in Lemma \ref{extendability}, we note that $(S^\vee/S)/H_S \subseteq M^\vee /M$. Therefore, since the discriminant of $S$ is $-2^{10}$ and we assumed $p \neq 2$, we have $H_S = (S^\vee/S)$.
By Lemma \ref{identification} $\tilde{g}^*|_S$ acts on $H_S$ as an isometry of order divisible by $17$ resp. $31$. Then, Lemma \ref{extendability} implies that $\tilde{g}^*|_K$ acts on $H_K$ as an isometry divisible by the same orders. On the other hand, ${\rm rk}(K) = 12$ and $K$ is negative definite, so that $\tilde{g}^*|_K$ has finite order and $\varphi(\ord(\tilde{g}^*|_K)) \leq 12$. This contradiction yields Claim $3$.
\vspace{3mm}

\noindent
\textbf{Claim 4:} If $\Char(k) \neq 2$, there is no element of order $3,5$ or $7$ in $\overline{\Aut}(X)$.

\indent \indent Suppose there was and let $g \in \Aut(X)$ with image $\overline{g} \in \overline{\Aut}(X)$ such that $\overline{g}$ has order $3,5$ or $7$. By Proposition \ref{orbits}, $N(2)$ is not divisible by $\ord(\overline{g})$ and hence there is some $W_{E_{10}}(2)$-orbit of non-degenerate canonical $U_{[2]}$-sequences which is preserved by $g$. In particular, there is some $h \in W_{E_{10}}(2) \subseteq \Aut(X)$ such that $g \circ h$ preserves a non-degenerate canonical $U_{[2]}$-sequence $(f_1,f_2)$. Replace $g$ by $g \circ h$ and note that they induce the same element in $\overline{\Aut}(X)$. Then, $g$ preserves the decomposition $E_{10} = \langle f_1,f_2 \rangle \oplus E_8$, hence it has finite order. Moreover, we can replace $g$ by some even power of $g$ without changing $\ord(\overline{g})$, so we can assume without loss of generality that $g$ has odd order. Now, we lift $f_1$ and $f_2$ to classes of half-fibers $F_1,F_1',F_2$ and $F_2'$, where $F_i' \in  |F_i + K_X|$. Since $g$ has odd order, it preserves $F_1,F_1',F_2$ and $F_2'$ and since $g$ fixes the four points of intersection of these four curves, it has at least two smooth fixed points on each of these curves. All half-fibers of elliptic fibrations on $X$ are isomorphic to irreducible plane cubics, because we assumed $X$ to be unnodal. Since $g$ has odd order and at least two smooth fixed points on the half-fibers, it follows from the known structure of the fixed loci of automorphisms of cubic curves (see e.g. \cite{Silverman}), that $g$ has at least three fixed points on each half-fiber. This implies that the action of $g$ on the linear systems $|2F_1|$ and $|2F_2|$ has at least three fixed points and is therefore trivial. But for a general member $F$ of the pencil $|2F_1|$, the set $F \cap (F_2 \cup F_2')$ consists of $4$ distinct points. Since $g$ has odd order and preserves $F_2$ and $F_2'$, the induced automorphism $g|_F$ has at least $4$ fixed points and is therefore trivial. This is a contradiction.

\vspace{3mm}
\noindent \textbf{Claim 5:} If $\Char(k) = 2$, then $\overline{\Aut}(X)$ is a $2$-group.

\indent \indent If $\Char(k) = 2$, then $X$ is ordinary by our convention and thus so is $Y$ by a result of R. Crew \cite[Theorem 2.7]{Crew}. Therefore, by a recent result of T. Srivastava \cite[Theorem 4.11]{Srivastava}, all automorphisms of $Y$ lift to the canonical lift $Y_\eta$ of $Y$. Hence, Claim $5$ follows from Claim $1$ and $4$ and the equivariance of specializations of line bundles.
\qed

\bigskip

\begin{Remark}
We remark that the case in Claim $3$ really occurs, that is, there are unnodal Enriques surfaces $X$ with supersingular K3 cover $Y$ (apart from the unnodal classical and supersingular Enriques surfaces in characteristic $2$). In fact, even Artin invariant $1$ is possible for $Y$, as the following example shows:

Over $\CC$, singular K3 surfaces are dense in the moduli space of K3 surfaces covering Enriques surfaces by the density of Noether--Lefschetz loci (see e.g. \cite[Theorem 1.1]{Oguiso}), while the set of K3 surfaces covering unnodal Enriques surfaces is open (see Proposition \ref{unnodalisopen}). Hence, there is a singular K3 surface $Y$ over $\CC$ with unnodal Enriques quotient $\pi: Y \to X$.
This $Y$ has complex multiplication and thus $Y$, $X$ and $\pi$ can be defined over a number field $K$. Spreading out, we obtain integral models $\mathcal{Y}$ and $\mathcal{X}$ of $Y$ and $X$ over the ring of integers $\mathcal{O}_K$ of $K$. Now, again by Proposition \ref{unnodalisopen} and since the singular loci of $\mathcal{Y}$ and $\mathcal{X}$ are closed, there is a closed and hence finite set $S \subseteq \Spec(\mathcal{O}_K)$ such that for all places $\mathfrak{p} \in \Spec(\mathcal{O}_K) \setminus S =: U$ the reductions $\mathcal{Y}_{\mathfrak p}$ and $\mathcal{X}_{\mathfrak p}$ modulo $\mathfrak{p}$ are smooth and $\mathcal{X}_{\mathfrak p}$ is unnodal. On the other hand, by a result of I. Shimada \cite[Theorem 1]{Shimada}, there is an infinite subset $V \subseteq \Spec(\mathcal{O}_K)$ such that for all $\mathfrak{p} \in V$, the reduction $\mathcal{Y}_{\mathfrak p}$ is supersingular with Artin invariant $1$. Therefore, for all of the infinitely many primes $\mathfrak{p} \in U \cap V$, the K3 surface $\mathcal{Y}_{\mathfrak p}$ is a supersingular K3 surface of Artin invariant $1$ covering the unnodal Enriques surface $\mathcal{X}_{\mathfrak p}$. Note that both $\mathcal{Y}_{\mathfrak p}$ and $\mathcal{X}_{\mathfrak p}$ are defined over $\overline{\FF}_p$.

For more information on Enriques quotients of supersingular K3 surfaces, we refer the reader to \cite{Jang3}.
\end{Remark}

Back to Enriques surfaces with extra automorphisms, we will now deduce that Theorem \ref{2group} implies that we can apply the techniques of Section \ref{section2}.

\begin{Corollary}\label{goodgroup}
Let $X$ be an unnodal Enriques surface. Then, there is an elliptic fibration $f: X \to \PP^1$ such that the subgroup $\Aut(f) \subseteq \Aut(X)$ surjects onto $\overline{\Aut}(X)$. 
\end{Corollary}

\prf
By Theorem \ref{2group}, $\overline{\Aut}(X)$ is a $2$-group and by Proposition \ref{orbits} the number of $W_{E_{10}}(2)$-orbits of effective, primitive, and isotropic vectors $F \in \Num(X)$ is odd, so that there is at least one such orbit which is preserved by the whole group $\overline{\Aut}(X)$. For an elliptic fibration $f$ defined by an element in this orbit, the group ${\rm Aut}(f) \cdot W_{E_{10}}(2)$ maps onto $\overline{\Aut}(X)$, hence so does ${\rm Aut}(f)$.
\qed
\bigskip

Now that we know that there is an elliptic fibration $f: X \to \PP^1$ inducing every extra automorphism of $X$, we can apply the results of Section \ref{section2} to this fibration. The following lemma shows that every automorphism of $f$ that is non-trivial modulo $W_{E_{10}}(2)$ acts non-trivially on $\PP^1$. This will allow us to relate $\overline{\Aut}(X)$ to $\overline{\Aut}_2$ of the Jacobian of $f$, which we have essentially computed in Section \ref{section3}.

\begin{Lemma}\label{actionpicard}
Let $f: X \to \PP^1$ be an elliptic fibration of an unnodal Enriques surface. Then,
$$
{\rm MW}(J(f)) \rtimes \ZZ/2\ZZ \cong \ker(p_{\PP^1}: \Aut(f) \to \Aut(\PP^1)) \subseteq W_{E_{10}}(2).
$$
\end{Lemma}

\prf
Let us first show the isomorphism. We use the notation of Section \ref{section2}. Let $C$ be the generic fiber of $f$ and $E$ its Jacobian. Then, by definition, ${\rm MW}(J(f)) = E(k(t))$. As in Section \ref{section2}, we have $\ker(p_{\PP^1}) = \Aut_{C/k(t)}(k(t))$ and an exact sequence
\begin{equation} \label{ses}
\xymatrix{
0 \ar[r] &  E(k(t)) \ar[r] & \Aut_{C/k(t)}(k(t)) \ar[r] & Q_C(k(t)).
}
\end{equation}
Recall that $i(C) = 2$ by Section \ref{EllfibronEnriques} and that $E(k(t))$ contains no $2$-torsion element because $J(f)$ has no reducible fibers (see \cite[p.111]{ShiodaSchuett}). Thus, Lemma \ref{QC} shows that $Q_C(k(t))$ is an elementary $2$-group. As is well-known, inversion with respect to the group structure on $C(\overline{k(t)})$ is the unique involution in $Q_C(\overline{k(t)})$. Hence, $Q_C(k(t)) \subseteq \ZZ/2\ZZ$.
To show that $\Aut_{C/k(t)}(k(t)) \cong {\rm MW}(J(f)) \rtimes \ZZ/2\ZZ$, it suffices to find an involution $g$ in $\Aut_{C/k(t)}(k(t))$, for then $g \not \in E(k(t))$ and thus $Q_C(k(t)) \cong \ZZ/2\ZZ$ and the Sequence (\ref{ses}) splits.

Let $F_1$ be a half-fiber of $f$. By \cite[Theorem 3.4.1]{CossecDolgachev} and Example \ref{BiellipticInvolution}, there is a half-fiber $F_2$ of an elliptic fibration different from $f$ such that $|2F_1 + 2F_2|$ induces a bielliptic involution $g$ and this involution preserves the fibers of $f$. Hence, $g \in \Aut_{C/k(t)}(k(t))$ is an involution. This proves ${\rm MW}(J(f)) \rtimes \ZZ/2\ZZ \cong \ker(p_{\PP^1})$.

Now, note that $g \in W_{E_{10}}(2)$ by Theorem \ref{Allcock} and ${\rm MW}(J(f))$ acts trivially on $\Num(X)/(i(C)\Num(X)) = \Num(X)/2\Num(X)$ by Corollary \ref{actionhorizontalpicard}. Thus, ${\rm MW}(J(f)) \subseteq W_{E_{10}}(2)$ by the definition of $W_{E_{10}}(2)$. This shows that $\ker(p_{\PP^1}) \subseteq W_{E_{10}}(2)$.
\qed

\bigskip

%

\smallskip

\begin{Corollary}[\textbf{$\mathrel{\widehat{=}}$ Theorem \ref{mainA}}] \label{CorollaryA}
Let $X$ be an unnodal Enriques surface. Then, $|\overline{\Aut}(X)| \in \{1,2,4\}$. If $\Char(k) = 2$, even $\overline{\Aut}(X) = \{1\}$ holds.
\end{Corollary}

\prf
By Corollary \ref{goodgroup}, we can choose an elliptic fibration $f: X \to \PP^1$ such that $\Aut(f)$ maps onto $\overline{\Aut}(X)$. Let $G := \varphi(\Aut(f)) \subseteq \Aut(J(f))$ be the image of the natural map from Corollary \ref{inducedactionpicard} and let $\overline{G} = p_{\PP^1}(G) \subseteq \Aut(\PP^1)$. Note that $\overline{G}$ is finite, since it permutes the points on $\PP^1$ lying under singular fibers and half-fibers of $f$ and there are finitely many but at least three such points by Remark \ref{atleastonesingfiber}. Applying Corollary \ref{inducedactionpicard}, Lemma \ref{actionpicard} and the assumption that $\Aut(f)$ maps onto $\overline{\Aut}(X)$, we see that 
$$
\overline{G} \cong p_{\PP^1}(\Aut(f)) \cong \Aut(f)/{\rm Ker}(p_{\PP^1}) \twoheadrightarrow \Aut(f)/(W_{E_{10}}(2) \cap \Aut(f)) \cong \overline{\Aut}(X).
$$
Since $\overline{\Aut}(X)$ is a $2$-group by Theorem \ref{2group}, a $2$-Sylow subgroup $\overline{G}_2$ of $\overline{G}$ still maps onto $\overline{\Aut}(X)$. 

If $\Char(k) \neq 2$, we can assume that the double fibers of $f$ lie over $0$ and $\infty$. Then, $\overline{G}_2$ is a $2$-group contained in $\overline{\Aut}_2(J(f),\sigma,\{0,\infty\})$, since $\overline{G} \cong G/(\ZZ/2\ZZ)$ by Lemma \ref{actionpicard} and Corollary \ref{inducedactionpicard}. Hence, we deduce that $
|\overline{\Aut}(X)|$ divides $|\overline{G}_2|$ which in turn divides $4$ by Corollary \ref{extrabound}.

If $\Char(k) = 2$, note that $\overline{G}_2$ fixes the point $p \in \PP^1$ lying under the unique double fiber of $f$. By the facts recalled in Section \ref{EllfibronEnriques}, this double fiber is of type $\I_0$ and ordinary or of type $\I_1$, hence the $j$-invariant of the Jacobian fibration $J(f)$ is not identically zero. Then, it is easy to see that there is a unique point $q \in \PP^1$ where the $j$-invariant of $J(f)$ vanishes (see also \cite{Lang1}). Thus, $\overline{G}_2$ fixes $p$ and $q$ and therefore $\overline{G}_2$ is a subgroup of $\GG_m$. Since $\overline{G}_2$ is a $2$-group and $\GG_m$ does not contain any element of order $2$, we have $\overline{G}_2 = \overline{\Aut}(X) = \{1\}$.
\qed

\section{Deformation theory and Moduli} \label{section5}
Throughout this section, we assume $\Char(k) \neq 2$. We will study deformations and families of unnodal Enriques surfaces and relate them to deformations of rational elliptic surfaces. Due to the work we have done in Section \ref{section3}, we have a complete understanding of the latter and this knowledge will allow us to prove Theorem \ref{mainB}.

\begin{Definition}
A \emph{family of (unnodal) Enriques surfaces} is a smooth and proper morphism  of algebraic spaces $\varphi: \mathcal{X} \to \mathcal{S}$ such that all geometric fibers of $\varphi$ are (unnodal) Enriques surfaces.
\end{Definition}

The following is certainly well-known (see e.g. \cite[Section 3]{CossecDolgachev2}), but we have not been able to find a reference in this generality.

\begin{Proposition}\label{unnodalisopen}
Let $\varphi: \mathcal{X} \to  \mathcal{S}$ be a family of Enriques surfaces. Let $U \subseteq |\mathcal{S}|$ be the subset of all points $s \in |\mathcal{S}|$ such that the geometric fiber $\mathcal{X}_{\overline{s}}$ is an unnodal Enriques surface. Then, $U$ is open in $|\mathcal{S}|$.
\end{Proposition}

\proof
Let $s \in U$. After replacing $\mathcal{S}$ by a sufficiently small fppf neighborhood of $s$, we can assume that $\mathcal{S}$ is an integral, affine scheme, that $\Pic_{ \mathcal{X} /\mathcal{S}}(\mathcal{S}) = \Pic(\mathcal{X})/\Pic(\mathcal{S})$ and that  $\Num_\mathcal{X} := \Pic_{\mathcal{X}/\mathcal{S}} / \Pic_{\mathcal{X}/\mathcal{S}}^\tau$ is the constant $\mathcal{S}$-group scheme whose fibers, together with the usual intersection pairing, are isomorphic to $E_{10}$ (see \cite[Corollary 4.3]{EkedahlHylShep}). Since $\Pic_{\mathcal{X}/\mathcal{S}}^\tau$ is an $\mathcal{S}$-group scheme of length $2$, any vector $\delta \in \Num_\mathcal{X}(\mathcal{S}) = E_{10}$ has at most two lifts $\delta_1,\delta_2$ to $\Pic_{\mathcal{X}/\mathcal{S}}(\mathcal{S})$. Choose line bundles $\mathcal{L}_1,\mathcal{L}_2 \in \Pic(\mathcal{X})$ representing $\delta_1,\delta_2$ and define
$$
T_{\delta} := {\rm Supp} (\varphi_* \mathcal{L}_1) \cup {\rm Supp} (\varphi_* \mathcal{L}_2),
$$
which is well-defined since $\Pic_{\mathcal{X}/\mathcal{S}}(\mathcal{S}) = \Pic(\mathcal{X})/\Pic(\mathcal{S})$ and closed in $\mathcal{S}$ since $\varphi$ is proper. Choose a line bundle $\mathcal{L} \in \Pic(\mathcal{X})$ of some degree $d$ which is relatively very ample over an open neighborhood $U'$ of $s$. By \cite[Main Theorem]{CossecDolgachev2}, we have
$$
U \cap U' = U' \setminus \left( \bigcup_{\delta \in I} T_\delta \right)
$$
where $I = \{\delta \in E_{10} \mid \delta^2 = -2 \text{ and } \delta.\mathcal{L} \leq d \}$. Note that $U \cap U'$ is an open neighborhood of $s$ in $U$ since the set $I$ is finite. Hence, $U$ is open.
\qed

\begin{Remark}
Note that the proof of Proposition \ref{unnodalisopen} works for arbitrary families of Enriques surfaces, not only for families over $k$ with $\Char(k) \neq 2$.
\end{Remark}

In Section \ref{section4} we have used the Jacobian fibration of an elliptic fibration of an unnodal Enriques surface $X$ to determine the possible orders of $\overline{\Aut}(X)$. The following lemma shows that it is even possible to relate the deformation theory of $X$ to the deformation theory of the Jacobian fibrations of elliptic fibrations on $X$.

\begin{Lemma} \label{complicated}
Let $\varphi: \mathcal{X} \to \mathcal{S}$ be a family of unnodal Enriques surfaces and assume that $\mathcal{S}$ is locally of finite type over $k$. Let $s \in \mathcal{S}(k)$ be a point with fiber $X := \mathcal{X}_s$ such that $\varphi$ is universal at $s$. Then, for every elliptic fibration $f: X \to \PP^1$ on $X$, there is an \'etale neighborhood $T$ of $s$ in $\mathcal{S}$ such that
\begin{enumerate}[(i)]
\item the elliptic fibration $f$ extends to a family of elliptic fibrations $\overline{f}: \mathcal{X}_T \to \PP^1_T$ over $T$,
\item there is a family of $2$-marked rational elliptic surfaces $(J(\overline{f}),\sigma,\{p_1,p_2\})$ over $T$ such that for all $t \in |T|$, the elliptic fibration $J(\overline{f})_{\overline t}$ together with the section $\sigma_{\overline t}$ is the Jacobian fibration of $\overline{f}_{\overline{t}}$, and the fibers of $\overline{f}_{\overline{t}}$ over $(p_1)_{\overline{t}}$ and $(p_2)_{\overline{t}}$ are the double fibers of $\overline{f}_{\overline{t}}$. In particular, all fibers of $(J(\overline{f}),\sigma,\{p_1,p_2\})$ satisfy $(\ast)$, and
\item the family $(J(\overline{f}),\sigma,\{p_1,p_2\})$ is universal at $s$.
\end{enumerate}
\end{Lemma}

\prf
Let us first prove $(i)$. By Grothendieck's existence theorem, it suffices to show that $f$ extends to a formal neighborhood of $\mathcal{X}_s$. Let $F_1,F_2$ be the half-fibers of $f$. Then, the morphism $f$ is determined by a choice of basis of $H^0(X,\mathcal{O}_X(2F_i))$, or equivalently, by a choice of $s_i \in H^0(X,\mathcal{O}_X(F_i))$, since the subspaces $H^0(X,\mathcal{O}_X(F_i)) \subseteq H^0(X,\mathcal{O}_X(2F_i))$ generate the space of global sections of $\mathcal{O}_X(2F_i)$.

Since $H^1(X,\mathcal{O}_X) = H^2(X,\mathcal{O}_X) = 0$, all line bundles on $X$ deform uniquely to every deformation of $X$ and since $F_1,F_2$ are the reduced curves underlying double fibers in $\Char(k) \neq 2$,  $\mathcal{O}_{F_i}(F_i)$ is a $2$-torsion line bundle on $F_i$ (see e.g. \cite[Chapter 7]{Badescu}). In particular, by tensoring the standard sequence for the divisor $F_i$ with $\mathcal{O}_X(F_i)$ and taking the long exact sequence in cohomology, we see that $H^1(X,\mathcal{O}_{X}(F_i)) = 0$ for $i = 1,2$. Thus, by \cite[Corollary 3.3.15]{Sernesi}, all global sections of $\mathcal{O}_{X}(F_i)$ lift to every deformation of $X$. Hence, also $f$ lifts to every deformation of $X$ and in particular to a formal neighborhood of $X_s$. Note that $T$ can be chosen in such a way that the family of base curves of the family of elliptic fibrations is isomorphic to $\PP^1_T$, since $\PP^1$ is infinitesimally rigid. 

As for $(ii)$, we can apply $(i)$ and assume without loss of generality that $f$ extends to a family of elliptic fibrations $\overline{f}: \mathcal{X} \to \PP^1_{\mathcal{S}}$. Moreover, since $H^2(X,T_X) = 0$, the deformation space ${\rm Def}_X \cong \hat{\mathcal{O}}_{S,s}$ is smooth over $k$, and thus we can further assume that $\mathcal{S}$ is a smooth and affine scheme over $k$. Let $p_1,p_2: \mathcal{S} \to \PP^1_S$ be the sections corresponding to the images of the families of double fibers of $\overline{f}$ in $\PP^1_{\mathcal{S}}$. Let $U \subseteq \PP^1_{\mathcal{S}}$ be the open complement of $p_1(\mathcal{S}),p_2(\mathcal{S})$ and let $\overline{f}_U: \mathcal{X}_U \to U$ be the restriction of $\overline{f}$ to $U$. Taking the relative compactified Jacobian (see \cite{altmankleiman}) of $\overline{f}_U$, we obtain a proper and flat morphism $J(\overline{f}_U): J(\mathcal{X}_U) \to U$ together with a section $\sigma_U: U \to J(\mathcal{X}_U)$, such that all geometric fibers of $J(\overline{f}_U)$ are integral curves of genus $1$ and such that $J(\mathcal{X}_U)$ is smooth over $\mathcal{S}$.

Let $\mathcal{L} = R^1\overline{f}_*\mathcal{O}_{\mathcal{X}}$. By \cite[Proposition 1.3]{LiuLorenziniRaynaud}, $\mathcal{L}|_U$ is canonically isomorphic to $R^1J(\overline{f}_U)_* \mathcal{O}_{J(\mathcal{X}_U)}$. Since $U$ is affine, it follows from \cite[Section 1]{Lang3}, that $J(\overline{f}_U)$ can be described as a closed subvariety of $\PP(\mathcal{O}_U \oplus \mathcal{L}^2|_U \oplus \mathcal{L}^{3}|_U)$ by a Weierstra{\ss} equation
$$
y^2 = x^3 + a_2x^2 + a_4x + a_6
$$
with $a_i \in H^0(U,\mathcal{L}^{-i}|_U)$. Now, note that for every $t \in |\mathcal{S}|$, the Jacobian $J(\overline{f}_{\overline t})$ of $\overline{f}_{\overline{t}}$ is the unique relatively minimal proper smooth model of $J(\overline{f}_U)_{\overline t}$ over $\PP^1_{\overline t}$ and hence the restrictions $(a_i)_{\overline{t}} \in H^0(U_{\overline{t}},\mathcal{L}^{-i}|_{U_{\overline{t}}})$ lift uniquely to $(b_i)_{\overline{t}} \in H^0(\PP^1_{\overline{t}},\mathcal{L}^{-i}|_{\PP^1_{\overline{t}}})$ and the equation
$$
y^2 = x^3 + b_2x^2 + b_4x + b_6
$$
defines $J(\overline{f}_{\overline{t}})$. Hence, to prove claim $(ii)$, it suffices to show that all these $(b_i)_{\overline t}$ glue together, i.e., that the $a_i$ are restrictions of sections $b_i \in H^0(\PP^1_{\mathcal S},\mathcal{L}^{-i})$.

To prove this, we consider a family of K3 covers $\mathcal{Y} \to \mathcal{X} \to \PP^1_{\mathcal{S}}$. Consider the Stein factorization $h$ of this composition. After possibly replacing $\mathcal{S}$ by an \'etale neighborhood of $s$, we can blow down $(-2)$-curves in the fibers of $h$ (see \cite[Section 3]{Artin2}). Thus, we can assume without loss of generality that $h: \mathcal{Y} \to \PP^1_\mathcal{S}$ has only irreducible fibers. Considering its relative compacified Jacobian, we get commutative diagrams
$$
\begin{minipage}[b]{0.3\linewidth}
\xymatrix{
\mathcal{Y} \ar[d]^{h} & \mathcal{X} \ar[d]^{\overline{f}} \\
\PP^1_{\mathcal{S}} \ar[r]^{j} & \PP^1_{\mathcal{S}}
}
\end{minipage}
\hspace{2cm}
\begin{minipage}[b]{0.3\linewidth}
\xymatrix{
J(\mathcal{Y}) \ar[d]^{J(h)} & J(\mathcal{X}_U) \ar[d]^{J(\overline{f}_U)} \\
\PP^1_{\mathcal{S}} \ar[r]^{j} & \PP^1_{\mathcal{S}}
}
\end{minipage}.
$$
The morphism $j$ is of degree $2$ and ramified exactly over $p_1$ and $p_2$ (see \cite[Proof of Lemma 2.17]{Martin}) and since taking Jacobians commutes with changing the base, the restriction of $J(h)$ to $j^{-1}(U)$
coincides with the base change of $J(\overline{f}_U)$ along $j$. Then, pullback along $j$ induces a map of long exact sequences for cohomology with support
$$
\xymatrix{
0 \ar[r] & H^0(\PP^1_{\mathcal S},j^{*}(\mathcal{L}^{i})) \ar[r] & H^0(j^{-1}(U),j^*(\mathcal{L}^{i})|_{j^{-1}(U)}) \ar^-o[r] & H^1_{ \{j^{-1}(p_1),j^{-1}(p_2)\}}(\PP^1_{\mathcal S},j^{*}\mathcal{L}^{i}) \ar[r] & \hdots \\
0 \ar[r] & H^0(\PP^1_{\mathcal S},\mathcal{L}^{i}) \ar[r]  \ar^{j^*}[u] & H^0(U,\mathcal{L}^{i}|_U) \ar^-o[r]  \ar^{j_U^*}[u] & H^1_{ \{p_1,p_2\}}(\PP^1_{\mathcal S},\mathcal{L}^{i}) \ar[r]  \ar^{j^*}[u]& \hdots
}.
$$
Since $j$ is a finite morphism of degree prime to $\Char(k)$, there is a non-degenerate trace map and therefore a splitting of $\mathcal{O}_{\PP^1_\mathcal{S}} \to j_*\mathcal{O}_{\PP^1_\mathcal{S}}$, hence all of the vertical arrows in the above diagram are injective. By the above diagram, the existence of sections $b_i$ restricting to $a_i$ is equivalent to $o(a_i) = 0$. However, we certainly have $o(j^*_U(a_i)) = 0$, since the base change of $J(\overline{f}_U)$ along $j$ extends to the family of elliptic surfaces $J(h)$ and the coefficients of a Weierstra{\ss} model of $J(h)$ give a lifting of $j^*_U(a_i)$ to $H^0(\PP^1_{\mathcal S},j^{*}(\mathcal{L}^{i}))$. Since $j^*$ is injective, we deduce that $o(a_i) = 0$. This proves claim $(ii)$.

To prove claim $(iii)$, we have to show that the morphism $\Phi: \Spf \hat{\mathcal{O}}_{\mathcal{S},s} \to {\rm Def}_{(J(\overline{f})_s,\sigma_s,\{(p_1)_s,(p_2)_s\})}$ induced by $(J(\overline{f}),\sigma,\{p_1,p_2\})$ is an isomorphism. Since $h^1(X,T_X) = 10$ and $h^2(X,T_X) = 0$, we have $\hat{\mathcal{O}}_{\mathcal{S},s} \cong k \llbracket t_1,\hdots,t_{10} \rrbracket$ and by Lemma \ref{Weierstrassuni} the same is true for the $k$-algebra that prorepresents the functor ${\rm Def}_{(J(\overline{f})_s,\sigma_s,\{(p_1)_s,(p_2)_s\})}$. Thus, it suffices to show that $\Phi$ is injective on tangent spaces, that is, that every morphism $\Spec(k[\epsilon]) \to \mathcal{S}$ mapping the closed point to $s$ such that the base change of $(J(\overline{f})_\epsilon,\sigma_\epsilon,\{(p_1)_\epsilon,(p_2)_\epsilon\})$ along this morphism is isomorphic to the trivial deformation is constant. Hence, let us choose a morphism $\Spec(k[\epsilon]) \to \mathcal{S}$ such that $(J(\overline{f})_\epsilon,\sigma_\epsilon,\{(p_1)_\epsilon,(p_2)_\epsilon\})$ is isomorphic to the trivial deformation. By construction, the deformation $(J(\overline{f})_\epsilon,\sigma_\epsilon,\{(p_1)_\epsilon,(p_2)_\epsilon\})$ comes from a deformation $\overline{f}_{\epsilon}: \mathcal{X}_{\epsilon} \to \PP^1_{\epsilon}$ of $\overline{f}_s: \mathcal{X}_s \to \PP^1$. Since $\varphi$ is universal at $s$, it suffices to show that $\mathcal{X}_{\epsilon}$ is trivial, or equivalently, that the induced deformation $\mathcal{Y}_{\epsilon}$ of the K3 cover $\mathcal{Y}_s$ of $\mathcal{X}_s$ is trivial.

As in the proof of claim $(ii)$, $\mathcal{Y}_{\epsilon}$ admits an elliptic fibration $h_{\epsilon}$ lifting the fibration $\overline{f}_{\epsilon}$. In particular, $h_{\epsilon}$ restricts to an elliptic fibration $h_s$ of $\mathcal{Y}_s$ and we assume without loss of generality that all fibers of $h_s$ are irreducible. The smooth locus of $h_s$, say $h_s^{\rm sm}: \mathcal{Y}_s^{\rm sm} \to \PP^1$, is a torsor under its Jacobian $J(h^{\sm}_s): \Pic^0_{\mathcal{Y}_s/\PP^1} \to \PP^1$ and hence corresponds to a class $\eta \in H^1_{\et}(\PP^1, \Pic^0_{\mathcal{Y}_s/\PP^1} )$. The smooth locus $h_{\epsilon}^{\rm sm}: \mathcal{Y}_{\epsilon}^{\rm sm} \to \PP^1_{\epsilon}$ of $h_{\epsilon}$ then corresponds to class $\eta_{\epsilon} \in H^1_{\et}(\PP^1, \Pic^0_{\mathcal{Y}_{\epsilon}/\PP^1_{\epsilon}})$ restricting to $\eta$. As explained for example in \cite[Section 4]{Partsch}, the restriction map $r$ sits inside a long exact sequence
\begin{equation} \label{restrictionsequence}
\xymatrix{
\hdots \ar[r]& H^1(\PP^1,{\rm Lie}(\Pic^0_{\mathcal{Y}_s/\PP^1})) \ar[r] & H^1_{\et}(\PP^1, \Pic^0_{\mathcal{Y}_{\epsilon}/\PP^1_{\epsilon}}) \ar[r]^r & H^1_{\et}(\PP^1, \Pic^0_{\mathcal{Y}_s/\PP^1} ) \ar[r] & \hdots
}
\end{equation}

Before we can make use of the above Sequence (\ref{restrictionsequence}), we have to show that $\Pic^0_{\mathcal{Y}_{\epsilon}/\PP^1_{\epsilon}}$ is isomorphic to the trivial deformation of $\Pic^0_{\mathcal{Y}_s/\PP^1}$. By assumption, $(J(\overline{f})_\epsilon,\sigma_\epsilon,\{(p_1)_\epsilon,(p_2)_\epsilon\})$ is trivial. Next, note that the deformation $\{(p_1)_\epsilon,(p_2)_\epsilon\}$ being trivial implies that the morphism $j_{\epsilon}: \PP^1_{\epsilon} \to \PP^1_{\epsilon}$ is the trivial deformation of $j_s$, since $j_{\epsilon}$ is uniquely determined by its branch locus. Therefore, the base change of $J(\overline{f})_\epsilon$ along $j_{\epsilon}$ is the trivial deformation of the base change of $J(\overline{f})_s$ along $j_s$. Since taking Jacobians commutes with changing the base, $\Pic^0_{\mathcal{Y}_{\epsilon}/\PP^1_{\epsilon}}$ is the smooth locus of the base change of $J(\overline{f})_\epsilon$ along $j_{\epsilon}$, and hence $\Pic^0_{\mathcal{Y}_{\epsilon}/\PP^1_{\epsilon}}$ is isomorphic to the trivial deformation of $\Pic^0_{\mathcal{Y}_s/\PP^1}$.

Thus, in Sequence (\ref{restrictionsequence}) there is a class $\eta_0 \in H^1_{\et}(\PP^1, \Pic^0_{\mathcal{Y}_{\epsilon}/\PP^1_{\epsilon}})$ corresponding to the trivial deformation of $\eta$. We have $r(\eta_0 - \eta_{\epsilon}) = (\eta - \eta) = 0$, hence $\eta_0 - \eta_{\epsilon}$ comes from $H^1(\PP^1,{\rm Lie}(\Pic^0_{\mathcal{Y}_s/\PP^1}))$ and in particular $\eta_0 - \eta_{\epsilon}$ is $p$-torsion with $p = \Char(k)$, since ${\rm Lie}(\Pic^0_{\mathcal{Y}_s/\PP^1})$ is a coherent sheaf on $\PP^1$. However, we claim that both $\eta_0$ and $\eta_{\epsilon}$ have order at most $2$. In fact, by Section \ref{EllfibronEnriques}, the family of elliptic fibrations $\overline{f}: \mathcal{X} \to \PP^1_{\mathcal{S}}$ admits a bisection after possibly shrinking $\mathcal{S}$. Hence, the family of elliptic fibrations $h: \mathcal{Y} \to \PP^1_\mathcal{S}$ of the K3 covers admits either a bisection or a section. In particular, since $\eta_{\epsilon}$ is the restriction of the cohomology class corresponding to the $\Pic^0_{\mathcal{Y}/\PP^1_{\mathcal{S}}}$-torsor given by the smooth locus of $h$, it has order at most $2$. The same is true for $\eta_0$ and hence $\eta_0 - \eta_{\epsilon}$ is $2$-torsion. Since $p \neq 2$, this shows $\eta_0 = \eta_{\epsilon}$, hence  $\mathcal{Y}_{\epsilon}^{\rm sm}$ is trivial and since the complement of $\mathcal{Y}_{s}^{\rm sm}$ in $\mathcal{Y}_{s}$ has depth at least $2$, $\mathcal{Y}_{\epsilon}$ is trivial. Therefore, also $\mathcal{X}_{\epsilon}$ is trivial, which is what we wanted to prove.
\qed

\bigskip

Finally, we can use this relation between deformations of unnodal Enriques surfaces and deformations of their Jacobians to determine the dimension of the moduli space of unnodal Enriques surfaces with extra automorphisms.

\begin{Corollary}[\textbf{ $\mathrel{\widehat{=}}$ Theorem \ref{mainB}}] \label{CorollaryB}
Let $\varphi: \mathcal{X} \to \mathcal{S}$ be a family of unnodal Enriques surfaces with $\mathcal{S}$ locally of finite type over $k$. Let $s \in \mathcal{S}(k)$ be a point such that $\varphi$ is universal at $s$. For $j \in \{2,4\}$, define the set
$$
\mathcal{Z}_j = \{t \in \mathcal{S}(k) ~|~ |\overline{\Aut}(X_{\overline t})| \geq j \} \subseteq |\mathcal{S}|.
$$
Then, there is an open neighborhood $s \in U \subseteq |\mathcal{S}|$ such that
\begin{itemize}
\item $U$ has dimension $10$,
\item $\overline{\mathcal{Z}_2} \cap U$ has dimension at most $5$, and
\item $\overline{\mathcal{Z}_4} \cap U$ has dimension at most $2$.
\end{itemize} 
In particular, a generic Enriques surface does not admit any extra automorphisms.
\end{Corollary}

\prf
Note that the statement is \'etale local on $\mathcal{S}$. Clearly, $U$ can be chosen to be $10$-dimensional since $h^1(\mathcal{X}_s,T_{\mathcal{X}_s}) = 10$ and $h^2(\mathcal{X}_s,T_{\mathcal{X}_s}) = 0$. By Proposition \ref{orbits}, there are $527$ elliptic fibrations on $\mathcal{X}_s$ up to the action of $W_{E_{10}}(2)$. Choose a set of representatives $\{f_i: \mathcal{X}_s \to \PP^1\}$, $i=1,\hdots,527$. By Lemma \ref{complicated}, we can assume without loss of generality that all $f_i$ extend to families of elliptic fibrations $\overline{f}_i: \mathcal{X} \to \PP^1_\mathcal{S}$ and since specialization of line bundles is equivariant with respect to the action of automorphism group, the $(\overline{f}_i)_{\overline{t}}$ form a full set of representatives for the $W_{E_{10}}(2)$-orbits of elliptic fibrations on $\mathcal{X}_{\overline t}$ for every $t \in |\mathcal{S}|$.

Denote the family of Jacobians of $\overline{f}_i$ as in Lemma \ref{complicated} (ii) by $(J(\overline{f}_i),\sigma_i,\{(p_1)_i,(p_2)_i\})$ and denote the closed subsets of $\mathcal{S}$ parametrizing points $t$ such that $\overline{\Aut}_2((J(\overline{f}_i),\sigma_i,\{(p_1)_i,(p_2)_i\})_{\overline t})$ contains a $2$-group of order at least $j$ by $\mathcal{Z}_j^i$. As in the proof of Corollary \ref{CorollaryA}, if $|\overline{\Aut}(X_{\overline t})| \geq j$, then there is an elliptic fibration $f_i$ such that $\overline{\Aut}_2((J(\overline{f}_i),\sigma_i,\{(p_1)_i,(p_2)_i\})_{\overline t})$ contains a $2$-group of order at least $j$. Hence, we have
$$
\mathcal{Z}_j  \subseteq \bigcup_{i = 1}^{527} \mathcal{Z}_j^i.
$$
Since the right hand side is closed of dimension $5$ if $j = 2$ (resp. $2$ if $j = 4$) by Theorem \ref{Weierstrassmoduli}, $\overline{\mathcal{Z}_j}$ is closed of dimension at most $5$ if $j = 2$ (resp. $2$ if $j =4$).
\qed

\begin{Remark}
As another immediate consequence of Lemma \ref{complicated}, we obtain that an arbitrary elliptic fibration $f: X \to \PP^1$ of a generic Enriques surface satisfies the following good properties, since the corresponding properties are satisfied for a generic $2$-marked rational elliptic surface, as can be checked explicitly using Weierstra{\ss} equations:
\begin{itemize}
\item All fibers of $f$ are irreducible.
\item $f$ has exactly $12$ reduced singular fibers, all of which are nodal cubics.
\item The half-fibers of $f$ are irreducible, smooth (and ordinary if $\Char(k) \neq 0$) elliptic curves.
\end{itemize}
\end{Remark}

\section{Examples} \label{section6}
We have seen in Theorem \ref{mainA} that $|\overline{\Aut}(X)| \in \{1,2,4\}$. The purpose of this section is to give examples showing that all possibilities actually occur. The following examples were found by H. Ohashi \cite{Ohashi} in the context of semi-symplectic automorphisms of Enriques surfaces over the complex numbers. In the following, we will explain how they give rise to Enriques surfaces with extra automorphisms also in positive characteristic. The most important observation is the following.

\begin{Lemma}\label{bicanonical}
Let $X$ be an unnodal Enriques surface and $g \in W_{E_{10}}(2) \subseteq \Aut(X)$. Then, $g$ acts trivially on $H^0(X,\omega_X^{\otimes 2})$.
\end{Lemma}

\prf
By Theorem \ref{Allcock}, it suffices to prove that a bielliptic involution acts trivially on $H^0(X,\omega_X^{\otimes 2})$. But this is clear since every involution acts trivially on $H^0(X,\omega_X^{\otimes 2})$.
\qed

\begin{Example}
Consider the double cover of $\PP^1 \times \PP^1$ defined by the affine equation
\begin{eqnarray*}
w^2 &=& z(A(y^4z^2 - z^2) + B(y^4z - z^3) + C(y^4-z^4) \\ && + D(y^3z^2 - yz^2) + E(y^3z - yz^3) + F(y^2z - y^2z^3))
\end{eqnarray*}
with parameters $A,B,C,D,E,F$. Ohashi shows that the automorphism
$$
g_1: (w,y,z) \mapsto \big( \frac{\sqrt{-1}w}{y^2z^3},\frac{1}{y},\frac{1}{z} \big)
$$
of order $4$ preserves the above equation and that, if $k = \CC$, the minimal model of a resolution of a generic member of the above family is an Enriques surface $X$ whose K3 cover has Picard rank $10$. In particular, $X$ is unnodal. Moreover, he shows that $g_1$ lifts to $X$ and acts as $-1$ on $H^0(X,\omega_X^{\otimes 2})$. In particular, by Lemma \ref{bicanonical}, $g_1$ induces an extra automorphism. Since the above family has $5$ moduli over $\CC$ (see \cite[Section 2]{Ohashi}), Corollary \ref{CorollaryB} shows that $\overline{\Aut}(X) = \ZZ/2\ZZ$.

In fact, the blow-ups and blow-downs needed to obtain an Enriques surface from a generic member in the above family can be realized over an open subset of $\Spec(R)$ for a finite extension $R$ of $\ZZ[A,B,C,D,E,F]$ (this follows from the explicit description of this process in \cite[Section 2]{Ohashi}). In particular, we obtain a family of Enriques surfaces over an open subset of $\Spec(R)$. Since the property of being unnodal is open in families of Enriques surfaces by Proposition \ref{unnodalisopen}, this yields examples of unnodal Enriques surfaces $X$ with $\overline{\Aut}(X) = \ZZ/2\ZZ$ in all but possibly  finitely many characteristics.
\end{Example}

\begin{Example}
Consider the double cover of $\PP^1 \times \PP^1$ defined by the affine equation
\begin{eqnarray*}
w^2 &=& z\big(A(y^4z^2 + \sqrt{-1}z^4 -z^2 - \sqrt{-1}y^4) + B(y^4z + \sqrt{-1} y^2z^3 - z^3 - \sqrt{-1}y^2z)\\  &&+ D(y^3z^2 + \sqrt{-1}yz^3 - yz^2 - \sqrt{-1}y^3z)\big)
\end{eqnarray*}
with parameters $A,B,D$. 
For a primitive $8$-th root of unity $\zeta_8$, the automorphism
$$
g_2: (w,y,z) \mapsto \left( \frac{\zeta_8 y^3w}{z^4},\frac{y}{z}, \frac{y^2}{z} \right)
$$
preserves the above equation and, if $k = \CC$, the minimal model of a resolution of a generic member of this family is an Enriques surfaces $X$ whose K3 cover has Picard rank $10$. In particular, $X$ is unnodal. Ohashi shows that the above family has $2$ moduli (see \cite[Section 2]{Ohashi}). Moreover, $g_2$ lifts to $X$ and acts as $\sqrt{-1}$ on $H^0(X,\omega_X^{\otimes 2})$. By Lemma \ref{bicanonical}, this means that $g_2$ induces an extra automorphism of order $4$, hence $\overline{\Aut}(X) = \ZZ/4\ZZ$ by Corollary \ref{CorollaryA}. As in the previous example, the above equation actually yields examples of unnodal Enriques surfaces $X$ with $\overline{\Aut}(X) = \ZZ/4\ZZ$ in all but possibly finitely many characteristics.
\end{Example}

We were unable to find examples of unnodal Enriques surfaces with $\overline{\Aut}(X) = (\ZZ/2\ZZ)^2$. The following proposition shows that it is impossible to find such $X$ with a K3 cover of Picard rank $10$, even for $k = \CC$. For this reason, even if we find a good candidate for an Enriques surface $X$ with $\overline{\Aut}(X) = (\ZZ/2\ZZ)^2$, it will not be as easy as in the above examples to prove that $X$ unnodal.

\begin{Proposition}
Let $X$ be an unnodal Enriques surface over $k = \CC$ such that $\overline{\Aut}(X)$ is not cyclic. Then, the Picard rank $\rho(Y)$ of its K3 cover $Y$ satisfies $\rho(Y) > 10$.
\end{Proposition}

\prf
Suppose $\rho(Y) = 10$. Let $T = \Pic(Y)^\perp \subset H^2(Y,\ZZ)$ be the transcendental lattice of $Y$. Then, $\overline{\Aut}(X)$ acts as a group of Hodge isometries on $T$, hence via a cyclic quotient by \cite[Corollary 3.4]{Huybrechts}. The assumption $\rho(Y) = 10$ implies $\pi^*(\Num(X)) = \Pic(Y)$. Hence, by Lemma \ref{identification}, the action of $\overline{\Aut}(X)$ on $(\pi^*(\Num(X)))^\vee/(\pi^*(\Num(X))) \cong T^\vee/T$ is faithful and therefore $\overline{\Aut}(X)$ is cyclic.
\qed

\bibliographystyle{alpha}
\bibliography{ExtraAut}
\end{document}